\documentclass[reqno,11pt]{amsart}
\usepackage[margin=1in]{geometry}

\usepackage{amsthm, amsmath, amssymb, bm}
\usepackage{tikz}
\usepackage{xcolor}

\usepackage{microtype} 

\usepackage[utf8]{inputenc}
\usepackage[T1]{fontenc}

\usepackage[textsize=small,backgroundcolor=orange!20]{todonotes}

\usepackage[hidelinks]{hyperref}
\usepackage{url}

\usepackage[noabbrev,capitalize]{cleveref}
\crefname{equation}{}{}

\usepackage{xcolor,graphics}

\numberwithin{equation}{section}

\newtheorem{theorem}{Theorem}[section]
\newtheorem{proposition}[theorem]{Proposition}
\newtheorem{lemma}[theorem]{Lemma}

\newtheorem{corollary}[theorem]{Corollary}
\newtheorem{conjecture}[theorem]{Conjecture}
\newtheorem*{question*}{Question} \Crefname{question}{Question}{Questions}

\theoremstyle{definition}
\newtheorem{definition}[theorem]{Definition}

\newtheorem{openproblem}[theorem]{Open problem}
\newtheorem{question}[theorem]{Question}

\theoremstyle{remark}
\newtheorem*{remark}{Remark}

\newcommand{\abs}[1]{\left\lvert#1\right\rvert}
\newcommand{\norm}[1]{\left\lVert#1\right\rVert}

\newcommand{\paren}[1]{\left( #1 \right)}
\newcommand{\sqb}[1]{\left[ #1 \right]}
\newcommand{\set}[1]{\left\{ #1 \right\}}

\newcommand{\avg}{\overline}

\newcommand{\cc}[4]{\left\langle \begin{matrix}
	#1 & #2 \\ #3 & #4
\end{matrix}\right\rangle}

\newcommand{\ccs}[4]{\left\langle \begin{smallmatrix}
	#1 & #2 \\ #3 & #4
\end{smallmatrix}\right\rangle}

\newcommand{\bx}{\bm{x}}
\newcommand{\T}{\intercal}
\newcommand{\om}{\ominus}

\newcommand{\EE}{\mathbb{E}}
\newcommand{\EEE}{\mathop{\mathbb{E}}}

\newcommand{\RR}{\mathbb{R}}

\newcommand{\NN}{\mathbb{N}}

\tikzstyle{P} = [draw, circle, black, fill, inner sep = 0pt, minimum width = 3pt]
\tikzstyle{every loop} = []
\newcommand{\tikzHind}{
  \begin{tikzpicture}[baseline, yshift=1pt]
    \path[use as bounding box] (-.15,-.1) rectangle (.6,.35);
    \draw (0.5,0) node[P] {} -- (0,0) node[P] {} edge[-,in = 45, out = 135, loop] ();
  \end{tikzpicture}
}

\newcommand{\tikzHtwoloops}{
  \begin{tikzpicture}[baseline,yshift=1pt]
    \path[use as bounding box] (-.15,-0.1) rectangle (0.55,.35);
    \draw (0,0) node[P] {} edge[-,in = 45, out = 135, loop] () 
         (.4,0) node[P] {} edge[-,in = 45, out = 135, loop] ();
  \end{tikzpicture}
}

\newcommand{\tikzHwr}{
  \begin{tikzpicture}[baseline,yshift=1pt]
    \path[use as bounding box] (-.15,-0.1) rectangle (0.95,.35);
    \draw (0,0) node[P] {} edge[-,in = 45, out = 135, loop] () 
      -- (.4,0) node[P] {} edge[-,in = 45, out = 135, loop] ()
      -- (.8,0) node[P] {} edge[-,in = 45, out = 135, loop] ();
  \end{tikzpicture}
}

\title{A reverse Sidorenko inequality}

\author[Sah]{Ashwin Sah}
\address{Massachusetts Institute of Technology, Cambridge, MA, USA}
\email{asah@mit.edu}

\author[Sawhney]{Mehtaab Sawhney}
\address{Massachusetts Institute of Technology, Cambridge, MA, USA}
\email{msawhney@mit.edu}

\author[Stoner]{David Stoner}
\address{Department of Matheamtics, Stanford University, Stanford, CA, USA}
\email{dwstoner@stanford.edu}

\author[Zhao]{Yufei Zhao}
\address{Department of Mathematics, Massachusetts Institute of Technology, Cambridge, MA, USA}
\email{yufeiz@mit.edu}
\thanks{YZ was supported by NSF Awards DMS-1362326 and DMS-1764176, and the MIT Solomon Buchsbaum Fund.}

\date{September 2018 (initial); January 2020 (revised)}

\begin{document}

\begin{abstract}
Let $H$ be a graph allowing loops as well as vertex and edge weights. We prove that, for every triangle-free graph $G$ without isolated vertices, the weighted number of graph homomorphisms $\hom(G, H)$ satisfies the inequality
\[
\hom(G, H ) \le \prod_{uv \in E(G)} \hom(K_{d_u,d_v}, H )^{1/(d_ud_v)},
\]
where $d_u$ denotes the degree of vertex $u$ in $G$.
In particular, one has 
\[
\hom(G, H )^{1/|E(G)|} \le \hom(K_{d,d}, H )^{1/d^2}
\]
for every $d$-regular triangle-free $G$. The triangle-free hypothesis on $G$ is best possible. More generally, we prove a graphical Brascamp--Lieb type inequality, where every edge of $G$ is assigned some two-variable function. These inequalities imply tight upper bounds on the partition function of various statistical models such as the Ising and Potts models, which includes independent sets and graph colorings.

For graph colorings, corresponding to $H = K_q$, we show that the triangle-free hypothesis on $G$ may be dropped; this is also valid if some of the vertices of $K_q$ are looped. A corollary is that among $d$-regular graphs, $G = K_{d,d}$ maximizes the quantity $c_q(G)^{1/|V(G)|}$ for every $q$ and $d$, where $c_q(G)$ counts proper $q$-colorings of $G$.

Finally, we show that if the edge-weight matrix of $H$ is positive semidefinite, then
\[
\hom(G, H) \le \prod_{v \in V(G)} \hom(K_{d_v+1}, H )^{1/(d_v+1)}.
\]
This implies that among $d$-regular graphs, $G = K_{d+1}$ maximizes $\hom(G, H)^{1/|V(G)|}$. For 2-spin Ising models, our results give a complete characterization of extremal graphs: complete bipartite graphs maximize the partition function of 2-spin antiferromagnetic models and cliques maximize the partition function of ferromagnetic models.

These results settle a number of conjectures by Galvin--Tetali, Galvin, and Cohen--Csikv\'ari--Perkins--Tetali, and provide an alternate proof to a conjecture by Kahn.
\end{abstract}

\maketitle

\section{Introduction} \label{sec:intro}

\subsection{Independent sets, colorings, and graph homomorphisms}

Consider the following extremal questions. Given a graph $G$, let $i(G)$ denote the number of its independent sets, $c_q(G)$ the number of its proper $q$-colorings\footnote{A \emph{proper $q$-coloring} of $G$ is an assignment of each vertex of $G$ to $[q] :=\{1, \dots, q\}$ so that no two adjacent vertices are assigned the same color (in particular, the colors are labeled).}, and $\hom(G, H)$ the number of its graph homomorphisms to $H$ (we allow $H$ to have loops,  and later, weights on its vertices and edges).\footnote{A \emph{graph homomorphism} from $G$ to $H$ is a map of vertices $\phi \colon V(G) \to V(H)$ such that $\phi(u)\phi(v)$ is an edge of $H$ whenever $uv$ is an edge of $G$.}

\begin{question} \label{q:ind}
	Fix $d$. Among $d$-regular graphs, which $G$ maximizes $i(G)^{1/|V(G)|}$?
\end{question}

\begin{question} \label{q:color}
	Fix $d$ and $q$. Among $d$-regular graphs, which $G$ maximizes $c_q(G)^{1/|V(G)|}$?
\end{question}

\begin{question} \label{q:hom}
	Fix $d$ and $H$. Among $d$-regular graphs, which $G$ maximizes $\hom(G,H)^{1/|V(G)|}$?
\end{question}

The third question encompasses the first two, as $i(G) = \hom(G, \tikzHind)$ and $c_q(G) = \hom(G, K_q)$. 

The exponential normalization is a natural choice. Indeed, replacing $G$ by a disjoint union of copies of itself does not change the quantity $\hom(G, H)^{1/|V(G)|}$, as $\hom(G_1 \sqcup G_2, H) = \hom(G_1, H) \hom(G_2, H)$, where $\sqcup$ denotes a disjoint union.

\cref{q:ind} was initially raised by Granville in 1988 in connection with the Cameron--Erd\H{o}s conjecture on the number of sum-free sets.  Alon~\cite{Alon91} and Kahn~\cite{Kahn01} conjectured that $G = K_{d,d}$ is the exact maximizer. Alon~\cite{Alon91} proved an asymptotic version as $d \to \infty$. Kahn~\cite{Kahn01} proved the exact version under the additional hypothesis that $G$ is bipartite using an elegant entropy argument that later became quite influential in the field. Zhao~\cite{Zhao10} later removed this bipartite assumption via a combinatorial argument. The results of Kahn~\cite{Kahn01} and Zhao~\cite{Zhao10} together answer \cref{q:ind}: the  maximizer is $K_{d,d}$ (unique up to taking disjoint unions of copies of itself).

Galvin and Tetali~\cite{GT04} initiated the study of \cref{q:color,q:hom} and extended Kahn's entropy method~\cite{Kahn01} to prove that, under the additional hypothesis that $G$ is bipartite, $G= K_{d,d}$ is also the maximizer for $\hom(G, H)^{1/|V(G)|}$. See Lubetzky and Zhao \cite[Section 6]{LZ15} for a different proof using H\"older/Brascamp--Lieb type inequalities. Can the bipartite hypothesis on $G$ also be dropped in this case? Not for all $H$: e.g., for $H = \tikzHtwoloops$, $G = K_{d+1}$ is the maximizer instead of $K_{d,d}$. Extending the technique for independent sets, Zhao~\cite{Zhao11} showed that the bipartite hypothesis can be dropped for certain classes of $H$, but the techniques failed for $H = K_q$, corresponding to colorings (\cref{q:color}). It remained a tantalizing conjecture to remove the bipartite hypothesis for colorings.

Recently, Davies, Jenssen, Perkins, and Roberts developed a novel technique called the ``occupancy method'' \cite{DJPR1}, which gave a new proof of the maximization problem for independent sets (\cref{q:ind}). Their method reduces the problem to a (potentially large) linear program. Applying their method, they gave a computer-assisted proof of the coloring conjecture (answering \cref{q:color}) for $d=3$ ~\cite{DJPR3}, later extended to $d=4$ by Davies~\cite{Dav}. The occupancy method was later extended to other applications concerning independent sets~\cite{DJPR2,PP}, as well as geometric applications concerning sphere packings~\cite{JJP1} and spherical codes~\cite{JJP2}. Despite its successes, the occupancy method has a number of drawbacks. Its progress on \cref{q:color} requires extremely rapidly growing computational resources for larger values of $d$, and furthermore, the method appears to be ill-suited for irregular graphs.

Here, we answer \cref{q:color} and show that $G=K_{d,d}$ is always the maximizer, thereby resolving the coloring conjecture.

\begin{theorem} \label{thm:color}
	Let $G$ be a $d$-regular graph and $q$ a positive integer. Then
	\[
		c_q(G)^{1/|V(G)|} \le c_q(K_{d,d})^{1/(2d)}.
	\]
\end{theorem}

We also prove a more general result for not necessarily regular graphs. It is analogous to our recent result~\cite{SSSZ1} for independent sets, which resolved Kahn's conjecture~\cite{Kahn01}. Here is a way to phrase the question. Instead of ranging over $d$-regular graphs, what if we range over all graphs with a fixed \emph{degree--degree distribution}, i.e., the distribution of the integer-pair $\{d_u, d_v\}$ over a uniform random edge $uv \in E(G)$, where $d_u$ is the degree of $u \in V(G)$? Kahn conjectured that, for independent sets, the maximizing $G$, conditioned on a fixed degree-degree distribution, remains a disjoint union of complete bipartite graphs of possibly different sizes. We recently proved Kahn's conjecture, resulting in the following theorem.

\begin{theorem}[\cite{SSSZ1}] \label{thm:irreg-ind}
	Let $G$ be a graph without isolated vertices. Let $d_v$ be the degree of vertex $v$ in $G$. Then
	\[
	i(G) \le \prod_{uv \in E(G)} i(K_{d_u,d_v})^{1/(d_u d_v)}.
	\]
\end{theorem}

Galvin~\cite{Gal06} conjectured (falsely) that \cref{thm:irreg-ind} could be extended to $\hom(\cdot, H)$ for every $H$ in place of $i(\cdot)$. Here we prove the extension for $H = K_q$, extending our \cref{thm:color} on the number of proper $q$-colorings to irregular graphs.

\begin{theorem} \label{thm:color-irreg}
	Let $G$ be a graph without isolated vertices, and $q$ a positive integer. Let $d_v$ be the degree of vertex $v$ in $G$. Then
	\[
	c_q(G) \le \prod_{uv \in E(G)} c_q(K_{d_u,d_v})^{1/(d_u d_v)}.
	\]
\end{theorem}

Let us state a more general version of \cref{thm:color-irreg} that interpolates between independent sets and proper colorings. Fix a finite set of \emph{colors} $\Omega$ as well as a subset $\Omega_\circ \subseteq \Omega$, called the \emph{looped colors}. A \emph{semiproper coloring} of $G$ is an assignment of each vertex of $G$ to $\Omega$ so that for every non-looped color (i.e., a color in $\Omega \setminus \Omega_\circ$), the set of vertices of $G$ of that color is an independent set. In other words, with $q = \abs{\Omega}$ and $\ell = \abs{\Omega_\circ}$, semiproper colorings correspond to homomorphisms from $G$ to $K_q^{\ell \circ}$, where $K_q^{\ell \circ}$ is the complete graph on $q$ vertices with exactly $\ell$ vertices looped. Proper colorings correspond to $\ell = 0$. Independent sets correspond to $(\ell, q) = (1,2)$. The following theorem interpolates between \cref{thm:irreg-ind,thm:color-irreg}.

\begin{theorem} \label{thm:semiproper}
	Let $G$ be a graph without isolated vertices, and $\ell \le q$ nonnegative integers. Let $d_v$ be the degree of vertex $v$ in $G$. Then
	\[
	\hom(G, K_q^{\ell \circ}) \le \prod_{uv \in E(G)} \hom(G, K_q^{\ell \circ})^{1/(d_u d_v)}.
	\]
\end{theorem}

Let us now move on to general graph homomorphisms. Here, \cref{q:hom} remains wide open. There has been a number of conjectures stated in the literature, though several of them have been falsified by counterexamples and then later revised~\cite{CCPT17,Gal06,Gal13,GT04,Ser18}. For example, it was first conjectured~\cite{GT04} that the maximizer is always $G = K_{d,d}$, and then later revised~\cite{Gal13} to $G \in \{K_{d+1},K_{d,d}\}$, though this was later shown false too~\cite{Ser18}. We do not even have a conjecture for what the set of possible maximizers $G$ is. It is even unknown whether the set of potentially maximizing $G$ is finite for each $d$. See the recent survey~\cite{Zhao17} for more discussion on this problem. 

It is natural to restrict $G$ in hope of a cleaner result. Cohen, Csikv\'ari, Perkins, and Tetali~\cite{CCPT17} conjectured that among \emph{triangle-free} graphs $G$, the maximizer is always $G = K_{d,d}$, extending the theorem of Galvin and Tetali~\cite{GT04} for bipartite $G$. We prove this conjecture.

\begin{theorem} \label{thm:K3-free-hom}
	Let $G$ be a triangle-free $d$-regular graph, and $H$ a graph allowing loops. Then
	\[
		\hom(G, H)^{1/|V(G)|} \le \hom(K_{d,d}, H)^{1/(2d)}.
	\]
\end{theorem}

We extend the result to irregular graphs and prove a corrected version of Galvin's conjecture~\cite{Gal06}.

\begin{theorem} \label{thm:K3-free-hom-irreg}
	Let $G$ be a triangle-free graph without isolated vertices, and $H$ a graph allowing loops. Then
	\[
		\hom(G, H) \le \prod_{uv \in E(G)} \hom(K_{d_u,d_v}, H)^{1/(d_ud_v)}.
	\]
\end{theorem}

\begin{remark}
\cref{thm:K3-free-hom-irreg} remains true even if $H$ has vertex and edge weights, so that $\hom(G, H)$ is interpreted as the partition function for a certain ``$H$-model'' on $G$ (e.g., the hard-core model generalizing independent sets, and the Potts model generalizing colorings). In fact, it follows by standard observations in graph limit theory~\cite{BCLS08,Lov} (namely, approximating a graphon by a sequence of $W$-random graphs) that the weighted and unweighted versions of \cref{thm:K3-free-hom-irreg} are actually equivalent.
\end{remark}

Furthermore, the triangle-free hypothesis is best possible in \cref{thm:K3-free-hom,thm:K3-free-hom-irreg}.

\begin{proposition}
	\label{prop:counterexample-triangle}
	For every graph $G$ with a triangle, there exists some graph $H$ so that the inequality in \cref{thm:K3-free-hom-irreg} is false.
\end{proposition}

The analogous minimization problem is also interesting and mysterious, though here we only mention a few known cases (see~\cite{Csi16ar}). For both independent sets ($H = \tikzHind$) \cite{CR14} and colorings ($H = K_q$) \cite[Lemma A.1]{BSVV08} (also see \cite[Theorem 8.3]{Zhao17}), the minimizer is $K_{d+1}$, whereas for the Widom--Rowlinson model ($H = \tikzHwr$), the ``minimizer'' is the infinite $d$-regular tree \cite{Csi16ar}.

\subsection{Graphons, norms, and reverse Sidorenko} \label{sec:graphon}

In the theory of graph limits~\cite{Lov}, a \emph{graphon} is a symmetric measurable function $W \colon \Omega \times \Omega \to [0,1]$ (symmetric means $W(x,y) = W(y,x)$), where $\Omega$ is some probability space. Define the \emph{$G$-density in $W$} by
\[
t(G, W) := \int_{\Omega^{V(G)}} \prod_{uv \in E(G)} W(x_u,x_v) \, d\bx_{V(G)} ,
\]
where $d\bx_{V(G)} := \prod_{v \in V(G)} dx_v$ is the product probability measure on $\Omega^{V(G)}$.

Every graph $H$ can be turned into a graphon $W_H \colon V(H) \times V(H) \to \{0,1\}$ by using the uniform probability measure on $V(H)$ and letting $W_H(x,y) = 1$ if $xy \in E(H)$, and $W_H(x,y)=0$ if $xy \notin E(H)$. Then $t(G, H) := t(G, W_H) = \hom(G, H)/|V(H)|^{|V(G)|}$ is the homomorphism density of $G$ to $H$. The graphon notation naturally allows us to consider edge and vertex weights on $H$.

\cref{thm:K3-free-hom,thm:K3-free-hom-irreg} are equivalent to the following graphon formulation: 
\[
	\label{eq:K3-free-W}
	t(G, W) \le \prod_{uv \in E(G)} t(K_{d_u,d_v}, W)^{1/(d_ud_v)},
\]
and in particular, for an $n$-vertex $d$-regular graph,
\[
	t(G, W) \le t(K_{d,d}, W)^{n/(2d)}.
\]

Let us write
\[
	\|W\|_G := \abs{t(G, W)}^{1/|E(G)|}.
\]
Despite the suggestive notation, $\norm{\cdot}_G$ is not always a norm. These quantities were first considered by Hatami~\cite{Hat10} in connection to Sidorenko's conjecture. See the recent work of Conlon and Lee~\cite{CL17} addressing the question of which graphs $G$ induce norms.

Our results above can now be written as
\[
t(G, W) \le \prod_{uv \in E(G)} \|W\|_{K_{d_u,d_v}},
\]
and, in particular, for $d$-regular graphs $G$,
\[
\norm{W}_G \le \norm{W}_{K_{d,d}}.
\]
In contrast, Sidorenko's conjecture says that for all bipartite graphs $G$, $t(G, W) \ge t(K_2, W)^{|E(G)|}$,
or equivalently $\|W\|_G \ge \|W\|_{K_2}$.
Sidorenko's conjecture~\cite{ES83,Sid93} has been proved for several families of graphs~\cite{BR65,CFS10,CKLL18,CL-blowup,Hat10,KLL16,LS,Sid93,SzeSidorenko}, though it remains open in general. The first open case of the conjecture is $G = K_{5,5}\setminus C_{10}$ (also known as the ``M\"obius strip'' graph, for it is the incidence graph for a simplicial complex model of the M\"obius strip viewed as gluing together five triangles). Sidorenko's conjecture proposes that $\|\cdot\|_{K_2}$ is a lower bound to $\|\cdot \|_G$, whereas our result proves an upper bound $\norm{\cdot}_{K_{d,d}}$ for triangle-free $d$-regular graphs $G$. It is for this reason that we give the name \emph{reverse Sidorenko inequality}.

\subsection{Graphical Brascamp--Lieb inequalities}
We prove a generalization of \cref{thm:K3-free-hom-irreg}, allowing possibly different two-variable functions on every edge of $G$. This generalization corresponds to  graph homomorphisms with list colorings, where every vertex of $G$ is assigned an ``allowable'' subset of vertices of $H$, and we  only consider homomorphisms assigning each vertex of $G$ to one of its allowable vertices of $H$. This generality is actually needed as a strong induction hypothesis for our proof.

From now on, $H$ will be a \emph{weighted graph}, which we define to be a symmetric measurable function $H \colon \Omega \times \Omega \to \RR_{\ge 0}$, where $\Omega$ is a measure space. We set\footnote{Such quantities are more commonly denoted $Z_H(G)$ for the partition function of a spin model with weights and interactions given by $H$. Here we prefer to extend $\hom(G,H)$ notation so as to be consistent with the case for simple graphs.}
\[
\hom(G, H) := \int_{\Omega^{V(G)}} \prod_{uv \in E(G)} H(x_u,x_v) \, d\bx_{V(G)}.
\]
Here $d\bx_{V(G)} = \prod_{v \in V(G)} dx_v$ and each $dx_v$ is the measure on $\Omega$, which is encoding vertex weights on $H$. Then \cref{thm:K3-free-hom,thm:K3-free-hom-irreg} hold for weighted graphs $H$ as well (see remark following \cref{thm:K3-free-hom-irreg}).

In \cref{sec:graphon}, in discussing graphons, it was important in the statement of Sidorenko's conjecture that $\Omega$ is a probability space, or else an extra normalizing factor is needed. In contrast, the inequalities that we prove in this paper are all scale-free in the sense that the measure of $\Omega$ does not have to be normalized.

The reader is welcome to think of $H$ as an edge-weighted graph (allowing loops) on a finite set of vertices $\Omega$ (the ``colors'') equipped with the counting measure. By a standard graph limit argument, this case is equivalent to the general result.

\begin{definition} \label{defn:graph-norm}For a two-variable function $f \colon \Omega_1 \times \Omega_2 \to \RR$, define 
\[
\|f\|_{K_{a,b}} := \Biggl\lvert \int_{\Omega_1^a \times \Omega_2^b} \prod_{\substack{1 \le i \le a \\ 1 \le j \le b}} f(x_i,y_j) \, dx_1 \dotsm d x_a dy_1 \dotsm dy_b \Biggr\rvert^{1/(ab)}.
\]
\end{definition}

This quantity (again, not always a norm) can be viewed as a bipartite analog of the graph ``norm'' earlier,  though here we do not require $f$ to be symmetric. We have monotonicity in that if $a\le c$ and $b \le d$, and $f$ is nonnegative-valued, then $\norm{f}_{K_{a,b}} \le \norm{f}_{K_{c,d}}$, which can be proved by two applications H\"older's inequality, as can be be seen after rewriting the above integral as $\int_{\Omega_1^a} \Bigl(\int_{\Omega_2} f(x_1,y)\cdots f(x_a,y) \,dy\Bigr)^b \, dx_1\cdots dx_a$.

The following theorem generalizes \cref{thm:K3-free-hom-irreg} upon taking the same $f_{uv} = H$ for all edges $uv$. 

\begin{theorem} \label{thm:graph-BL}
	Let $G = (V,E)$ be a triangle-free graph. Let $\Omega_v$ be a measure space for each vertex $v \in V$. For each edge $uv \in E$, let $f_{uv} \colon \Omega_u \times \Omega_v \to \RR_{\ge 0}$ be a measurable function, labeled so that $f_{uv}(x_u,x_v)=f_{vu}(x_v,x_u)$. We have
	\[
	\int_{\Omega_V} \prod_{uv\in E} f_{uv}(x_u,x_v) \, d\bm x_V
	\le \prod_{uv\in E} \|f_{uv}\|_{K_{d_v,d_u}},
	\]
	where $\Omega_V := \prod_{v \in G} \Omega_v$ and $d\bm x_V := \prod_{v\in V} dx_v$, and $d_u$ is the degree of $u$ in $G$.
\end{theorem}

\begin{remark}
We have equality if (1) $G$ is a disjoint union of complete bipartite graphs, or (2) if there are functions $g_v \colon \Omega_v \to \RR_{\ge 0}$ such that $f_{uv}(x,y) = g_u(x)g_v(y)$ for every $uv \in E$.

By \cref{prop:counterexample-triangle}, the triangle-free hypothesis cannot be weakened.
\end{remark}

For semiproper list colorings, \cref{thm:graph-BL} holds without the triangle-free hypothesis, generalizing \cref{thm:semiproper}. See \cref{sec:color} for the statement and proof.

\cref{thm:graph-BL} can be viewed as a graphical analog of the Brascamp--Lieb inequalities~\cite{BL76,Lieb90}, which have the form $\int f_1(B_1 \bm x) \cdots f_k(B_k \bm x) d\bm x \lesssim \|f_1\|_{L^{p_1}} \cdots  \|f_k\|_{L^{p_k}}$, where the $B_i$'s are linear maps. The Brascamp--Lieb inequalities generalize classical inequalities such as H\"older's inequality and the Loomis--Whitney inequality, and have far reaching applications. Our inequality bounds a certain graphical integral in terms of graphical norm-like quantities that are in general weaker than $L^p$ norms. It may be possible that these graphical Brascamp--Lieb inequalities have a rich theory yet to be uncovered, e.g., extensions to more general setups such as hypergraphs and simplicial complexes, allowing greater flexibility in the combinatorial form of the integral on the left-hand side of the inequality.

\subsection{Positive semidefinite models are clique-maximizing} 

We have stated various results affirming that $G = K_{d,d}$ maximizes $\hom(G, H)^{1/|V(G)|}$ under various circumstances, e.g., among triangle-free $G$, or if $H$ is a (partially looped) complete graph corresponding to (semi)proper colorings. However, as remarked following \cref{thm:semiproper}, $K_{d,d}$ is not always the correct answer to \cref{q:hom}, and the general question remains very much open.

Given a weighted graph $H \colon \Omega \times \Omega \to \RR_{\ge 0}$, we say that $H$ is \emph{biclique-maximizing} if it satisfies, for all graphs $G$ without isolated vertices,
\begin{equation} \label{eq:biclique-maximizing}
	\hom(G, H) \le \prod_{uv \in E(G)} \hom(K_{d_u,d_v}, H)^{1/(d_ud_v)},
\end{equation}
where, as usual, $d_v$ denotes the degree of $v$ in $G$. We say that $H$ is \emph{clique-maximizing} if it satisfies, for all graphs $G$, 
\begin{equation} \label{eq:clique-maximizing}
	\hom(G, H) \le \prod_{v \in V(G)} \hom(K_{d_v+1}, H)^{1/(d_v+1)}.
\end{equation}
\cref{thm:semiproper} says that the partially looped complete graphs $K_q^{\ell \circ}$ are biclique-maximizing. On the other hand, it is not hard to check that a disjoint union of loops is clique-maximizing. It is known that there are graphs $H$ that are neither biclique-maximizing nor clique-maximizing, even among $d$-regular graphs $G$ (it is unknown which $G$ achieves the maximum for such $H$) \cite{Ser18}.

It was shown~\cite{Zhao11} that certain graphs $H$ satisfy a ``bipartite swapping trick'' inequality: $\hom(G,H)^2 \le \hom(G \times K_2, H)$ for all graphs $G$ (extended to a larger class of $H$ in \cite{Ser18}; also see \cref{sec:2spin} for a weighted variant). Here $G \times K_2$ is the graph with vertex set $V(G) \times \{0,1\}$ and an edge between $(v,i)$ and $(u,1-i)$ for every $uv \in E(G)$ and $i \in \{0,1\}$. Every $H$ satisfying this inequality is biclique-maximizing, since we can apply \cref{thm:graph-BL} to upper bound $\hom(G \times K_2, H)$ as $G \times K_2$ is bipartite and hence triangle-free.

In \cite{CCPT17,CPT17,Ser18}, it was shown that the Widom--Rowlinson model ($H = \tikzHwr$) satisfies \eqref{eq:clique-maximizing} for $d$-regular graphs $G$ (this was the first and essentially only such non-trivial case that was known). However, it turns out that $H = \tikzHwr$ is actually not clique-maximizing among irregular graphs (a counterexample is $G = K_{1, 4}$, as $113 >  7^{4/2} 63^{1/5}$). This case is interesting as there is a different maximization behavior between regular and irregular $G$. 

\begin{openproblem} \label{open:determine}
Determine all biclique-maximizing graphs $H$ and all clique-maximizing graphs $H$, in each case, for $d$-regular $G$ as well as for all $G$.
\end{openproblem}

We say that a weighted graph (also called a \emph{model}) $H \colon \Omega \times \Omega \to \RR_{\ge 0}$ is \emph{positive semidefinite} or \emph{ferromagnetic} if the corresponding function is  positive semidefinite (equivalently, the matrix $(H(x_i, x_j))_{i,j\in[n]}$ is positive semidefinite for every $x_1, \dots, x_n \in \Omega$). We say that $H$ is \emph{antiferromagnetic} if all eigenvalues (counting multiplicities) other than the top one are nonpositive. These definitions were taken from \cite{GSV14}.

For example, a disjoint union of loops is ferromagnetic, whereas $K_q^{\ell \circ}$ is antiferromagnetic. For 2-spin models, i.e., $\Omega = \{0,1\}$ allowing vertex weights, $H$ is ferromagnetic if $H(0,0)H(1,1) \ge H(0,1)^2$, and antiferromagnetic if $H(0,0)H(1,1) \le H(0,1)^2$. 

We prove the following result. See \cref{thm:clique-max-list} for a list coloring type generalization.

\begin{theorem} \label{thm:clique-max}
	Every ferromagnetic (i.e., positive semidefinite) model is clique-maximizing.
\end{theorem}

We conjecture that the converse holds as well.

Every 2-spin model is either ferromagnetic or antiferromagnetic depending on the sign of the determinant of its $2 \times 2$ edge-weight matrix, though this is false for $k$-spin models for $k > 2$.
As a corollary, we completely characterize all 2-spin models, generalizing independent sets. See \cref{sec:2spin} for the antiferromagnetic part of the proof, which follows from the bipartite swapping trick~\cite{Zhao10,Zhao11} and \cref{thm:graph-BL}.

\begin{corollary} \label{cor:2spin}
	A 2-spin model is biclique-maximizing if it is antiferromagnetic and clique-maximizing if it is ferromagnetic.
\end{corollary}

We close with a conjecture generalizing \cref{thm:semiproper}.

\begin{conjecture} \label{conj:antiferro}
	Every antiferromagnetic model is biclique-maximizing.
\end{conjecture}

The converse of \cref{conj:antiferro} is false. The looped graph $H$ given by the adjacency matrix $\left(\begin{smallmatrix} 1 & 1 & 1 \\ 1 & 1 & 0 \\ 1 & 0 & 0 \end{smallmatrix}\right)$ satisfies the bipartite swapping trick from \cite{Zhao11} and hence is biclique-maximizing (see comment two paragraphs before \cref{open:determine}) but $H$ is not antiferromagnetic. It remains wide open to classify all biclique-maximizing $H$.

\cref{thm:semiproper} establishes \cref{conj:antiferro} for $K_q^{\ell \circ}$. Though, even the following extension remains just out of reach of our current methods: $H \colon \Omega \times \Omega \to [0,1]$ with $H(x,y) = 1$ if $x \ne y$ and $0\le H(x,x) \le 1$ arbitrary. This is a generalization of the antiferromagnetic Potts model. The usual Potts model has additionally the same diagonal values $H(x,x) = \beta \in [0,1]$ for all $x \in \Omega$, and for these $H$, the conjecture has been verified for $3$-regular~\cite{DJPR3} and $4$-regular~\cite{Dav} graphs $G$ via the occupancy method with computer assistance.

\subsection{Relation to previous work} This work builds on our earlier work~\cite{SSSZ1} proving Kahn's conjecture on independent sets, \cref{thm:irreg-ind}, but requires several significantly new ideas. Our proof of \cref{thm:graph-BL} in \cref{sec:triangle-free} actually gives a new and more streamlined proof of \cref{thm:irreg-ind}. The new proof is significantly shorter, and it replaces a number of fairly technical inequality verifications in \cite{SSSZ1} (often involving checking repeated derivatives) by more conceptual inequalities primarily relying on H\"older's inequality and log-convexity considerations. In \cite{SSSZ1}, as in the earlier \cite{GZ11}, we relied on the recurrence $i(G) = i(G - v) + i(G - v - N(v))$ for the number independent sets, but such a relation is unavailable for colorings. Assigning a color to a vertex restricts the colors available to the neighborhoods, so it is natural to study the problem in the greater generality of list colorings and state a stronger induction hypothesis. By considering the effect of fixing a color on a vertex and carefully bounding contributions from far away vertices, we reduce the problem to more ``local'' inequalities. \cref{sec:toy} of the paper discusses the general reduction to local inequalities in greater detail.

\subsection*{Organization} In \cref{sec:toy} we give a toy calculation illustrating some proof ideas. In \cref{sec:triangle-free}, we prove \cref{thm:graph-BL}, the graphical Brascamp--Lieb inequality, and hence \cref{thm:K3-free-hom,thm:K3-free-hom-irreg}. In \cref{sec:color}, we prove \cref{thm:semiproper} concerning semiproper colorings, and hence \cref{thm:color,thm:color-irreg}. In \cref{sec:clique}, we prove \cref{thm:clique-max} showing that ferromagnetic (i.e., positive semidefinite) models are clique-maximizing.

\tikzstyle{rect}=[rectangle,draw=black,fill=white]
\tikzstyle{rectE}=[rectangle,draw=white,fill=white]
\tikzstyle{rectY}=[rectangle,draw=black,fill=yellow]

\newcommand{\wcirc}{\tikz\draw[draw=none,fill=none] (0,0) circle (.5ex);}
\newcommand{\bcirc}{\tikz\draw[white!40!blue,fill=white!40!blue] (0,0) circle (.5ex);}
\newcommand{\gcirc}{\tikz\draw[black!30!green,fill=black!30!green] (0,0) circle (.5ex);}
\newcommand{\rcirc}{\tikz\draw[red,fill=red] (0,0) circle (.5ex);}

\newcommand{\lc}[3]
{
 \ifnum#1=0
  \wcirc{}
 \else
  \rcirc{}
 \fi
 \ifnum#2=0
  \wcirc{}
 \else
  \gcirc{}
 \fi
 \ifnum#3=0
  \wcirc{}
 \else
  \bcirc{}
 \fi
}

\newcommand{\wcircg}{\tikz\draw[draw=none,fill=none] (0,0) circle (.5ex);}
\newcommand{\bcircg}{\tikz\draw[white!80!blue,fill=white!80!blue] (0,0) circle (.5ex);}
\newcommand{\gcircg}{\tikz\draw[white!50!green,fill=white!50!green] (0,0) circle (.5ex);}
\newcommand{\rcircg}{\tikz\draw[white!70!red,fill=white!70!red] (0,0) circle (.5ex);}

\newcommand{\lcg}[3]
{
 \ifnum#1=0
  \wcircg{}
 \else
  \rcircg{}
 \fi
 \ifnum#2=0
  \wcircg{}
 \else
  \gcircg{}
 \fi
 \ifnum#3=0
  \wcircg{}
 \else
  \bcircg{}
 \fi
}
\tikzstyle{vcenter} = [baseline={([yshift=0ex]current bounding box.center)}]
\newcommand{\EQ}{
\begin{tikzpicture}[vcenter]
\node (A) at (0,1.4) {};
\node (B) at (0,0){{\huge \textbf{=}}};
\node (C) at (0,-1.4){};
\end{tikzpicture}
}

\newcommand{\LE}{
\begin{tikzpicture}[vcenter]
\node (A) at (0,1.4) {};
\node (B) at (0,0){\scalebox{1.5}{$\leq$}};
\node (C) at (0,-1.4){};
\end{tikzpicture}
}

\newcommand{\PLUS}{
\begin{tikzpicture}[vcenter]
\node (A) at (0,1.4) {};
\node (B) at (0,0){{\huge \textbf{+}}};
\node (C) at (0,-1.4){};
\end{tikzpicture}
}

\newcommand{\setupCoords}
{
\coordinate (v1) at (0,0);
\coordinate (v2) at (0.8,1.4);
\coordinate (v3) at (2.4,1.4);
\coordinate (v4) at (3.2,0);
\coordinate (v5) at (2.4,-1.4);
\coordinate (v6) at (0.8,-1.4);
\coordinate (v2A) at (0.6, 1.5);
\coordinate (v2B) at (0.4, 1.5);
\coordinate (v2C) at (0.4, 1.3);
\coordinate (v2D) at (0.6, 1.3);
\coordinate (v5A) at (0.6, -1.5);
\coordinate (v5B) at (0.4, -1.5);
\coordinate (v5C) at (0.4, -1.3);
\coordinate (v5D) at (0.6, -1.3);
\coordinate (v3A) at (0.97, 1.5);
\coordinate (v3B) at (1.17, 1.5);
\coordinate (v3C) at (1.13, 1.3);
\coordinate (v3D) at (0.97, 1.3);
\coordinate (v4A) at (1, -1.5);
\coordinate (v4B) at (1.2, -1.5);
\coordinate (v4C) at (1.2, -1.3);
\coordinate (v4D) at (1, -1.3);
}

\newcommand{\moreCoords}
{
\coordinate (c1) at (1.2, 1.8);
\coordinate (c2) at (4.6, 1.8);
\coordinate (c4) at (1.2, -1.8);
\coordinate (c3) at (4.6, -1.8);
\coordinate (c10) at (1.3, 1.8);
\coordinate (c20) at (4.5, 1.8);
\coordinate (c40) at (1.3, -1.8);
\coordinate (c30) at (4.5, -1.8);
\coordinate (v11) at (0,0.4);
\coordinate (v13) at (0,2.0);
\coordinate (v12) at (-0.8,1.2);
\coordinate (v14) at (0.8,1.2);
\coordinate (v21) at (2,2.0);
\coordinate (v23) at (3.6,2.0);
\coordinate (v22) at (2,1.2);
\coordinate (v24) at (3.6,1.2);
\coordinate (v25) at (3.6,1.6);
\coordinate (v21C) at (2,1.6);
\coordinate (v23C) at (3.6,1.6);
\coordinate (v22C) at (2,0.8);
\coordinate (v24C) at (3.6,0.8);
\coordinate (v21B) at (2,1.2);
\coordinate (v22B) at (2,0.4);
\coordinate (v25B) at (3.6,0.8);
\coordinate (v31) at (4.8,1.2);
\coordinate (v33) at (5.6,0.4);
\coordinate (v32) at (5.6,2.0);
\coordinate (v34) at (6.4,1.2);
\coordinate (p1) at (0.8,1.8);
\coordinate (p2) at (4.0,2.6);
\coordinate (p22) at (4.0,-2.6);
\coordinate (p2C) at (4.0,2.2);
\coordinate (p2CC) at (4.0,-2.2);
\coordinate (p2CCC) at (4.0,0.2);
\coordinate (p3) at (6.4,1.8);
\coordinate (p2B) at (4.8,2);
\coordinate (p2BB) at (4.8,-2);
\coordinate (p4) at (0.8,-0.6);
\coordinate (p5) at (4.0,-0.6);
\coordinate (p5C) at (4.0,-0.2);
\coordinate (p5CC) at (4.0,-2.2);
\coordinate (p6) at (6.4,-0.6);
\coordinate (v41) at (0,-0.4);
\coordinate (v43) at (0,-2.0);
\coordinate (v42) at (-0.8,-1.2);
\coordinate (v44) at (0.8,-1.2);
\coordinate (v51) at (2,-2.0);
\coordinate (v53) at (3.6,-2.0);
\coordinate (v52) at (2,-1.2);
\coordinate (v54) at (3.6,-1.2);
\coordinate (v51C) at (2,-1.6);
\coordinate (v53C) at (3.6,-1.6);
\coordinate (v52C) at (2,-0.8);
\coordinate (v54C) at (3.6,-0.8);
\coordinate (v51D) at (2,-1.6);
\coordinate (v53D) at (0.4,-1.2);
\coordinate (v52D) at (2,-0.8);
\coordinate (v54D) at (3.6,-1.2);
\coordinate (v55) at (3.6,-1.6);
\coordinate (v51B) at (2,-1.2);
\coordinate (v52B) at (2,-0.4);
\coordinate (v55B) at (3.6,-0.8);
\coordinate (v61) at (4.8,-1.2);
\coordinate (v63) at (5.6,-0.4);
\coordinate (v62) at (5.6,-2.0);
\coordinate (v64) at (6.4,-1.2);
\coordinate (x1) at (6.4,-2.0);
\coordinate (x2) at (4.8, 2.0);
\coordinate (x3) at (4.8,-2.0);
\coordinate (x4) at (6.4, 2.0);

}

\section{A toy calculation} \label{sec:toy}

In this section we sketch a toy calculation demonstrating our general inductive strategy on $G = C_6$ with target graph $K_3$, i.e., 3-list-coloring a 6-cycle. This is the dessert before the dinner, as the actual proof involves more difficult steps not shown here.

The inequality that we would like to prove is illustrated by the following diagram. This is a special case of \cref{thm:graph-BL} for list coloring, i.e, $f(x,y) = 1_{x \neq y}$ in \cref{thm:graph-BL} or $H = K_q$ in the graph homomorphism setup. See \cref{thm:semiproper-list} for a statement of the list coloring inequality. This is an example of the strong induction hypothesis for upper bounding the number of list colorings, and we will apply induction on the number of vertices of $G$.
\begin{equation}\label{eq:toy-c6}
\scalebox{0.7}{
\begin{tikzpicture}[vcenter]
\setupCoords{};
\draw[line width=1] (v1)--(v2)--(v3)--(v4)--(v5)--(v6)--(v1);
\node[rect] at (v1) {\lc101};
\node[rect] at (v2) {\lc110};
\node[rect] at (v3) {\lc011};
\node[rect] at (v4) {\lc111};
\node[rect] at (v5) {\lc101};
\node[rect] at (v6) {\lc111};
\end{tikzpicture}
\quad\LE\quad
\begin{tikzpicture}[vcenter]
\moreCoords{};
\draw[line width=1] (v11)--(v14)--(v12)--(v13)--(v11);
\draw[line width=1] (v21)--(v24)--(v22)--(v23)--(v21);
\draw[line width=1] (v31)--(v34)--(v32)--(v33)--(v31);
\draw[line width=1] (v41)--(v44)--(v42)--(v43)--(v41);
\draw[line width=1] (v51)--(v54)--(v52)--(v53)--(v51);
\draw[line width=1] (v61)--(v64)--(v62)--(v63)--(v61);
\node[rect] at (v11) {\lc101};
\node[rect] at (v12) {\lc101};
\node[rect] at (v13) {\lc110};
\node[rect] at (v14) {\lc110};
\node[rect] at (v21) {\lc110};
\node[rect] at (v22) {\lc110};
\node[rect] at (v23) {\lc011};
\node[rect] at (v24) {\lc011};
\node[rect] at (v31) {\lc011};
\node[rect] at (v32) {\lc011};
\node[rect] at (v33) {\lc111};
\node[rect] at (v34) {\lc111};
\node[rect] at (v41) {\lc101};
\node[rect] at (v42) {\lc101};
\node[rect] at (v43) {\lc111};
\node[rect] at (v44) {\lc111};
\node[rect] at (v51) {\lc111};
\node[rect] at (v52) {\lc111};
\node[rect] at (v53) {\lc101};
\node[rect] at (v54) {\lc101};
\node[rect] at (v61) {\lc101};
\node[rect] at (v62) {\lc101};
\node[rect] at (v63) {\lc111};
\node[rect] at (v64) {\lc111};
\node at (p1) {$\frac{1}{4}$};
\node at (p2) {$\frac{1}{4}$};
\node at (p3) {$\frac{1}{4}$};
\node at (p4) {$\frac{1}{4}$};
\node at (p5) {$\frac{1}{4}$};
\node at (p6) {$\frac{1}{4}$};
\node at (p22) {$\color{white} \frac{1}{4}$};
\end{tikzpicture}
}
\end{equation}

Let us explain the meaning of the above diagram. On the left-hand side, the figure should be interpreted as the number of valid list colorings of the 6-cycle where each vertex of the 6-cycle is assigned one of its listed colors, such that no adjacent vertices receive the same color. On the right hand side, we have a product of six quantities, each being the number of list colorings of a 4-cycle (with different color lists for each 4-cycle) raised to the power $1/4$.

To prove \eqref{eq:toy-c6}, we begin by selecting the color of the left-most vertex of the 6-cycle, which gives the following:
\begin{equation}\label{eq:toy-split-v0}
\scalebox{0.8}{
\begin{tikzpicture}[vcenter]
\setupCoords{};
\draw[line width=1] (v1)--(v2)--(v3)--(v4)--(v5)--(v6)--(v1);
\node[rectY] at (v1) {\lc101};
\node[rect] at (v2) {\lc110};
\node[rect] at (v3) {\lc011};
\node[rect] at (v4) {\lc111};
\node[rect] at (v5) {\lc101};
\node[rect] at (v6) {\lc111};
\end{tikzpicture}}
\quad =\quad
\scalebox{0.8}{\begin{tikzpicture}[vcenter]
\setupCoords{};
\draw[line width=1] (v2)--(v3)--(v4)--(v5)--(v6);
\node[rect] at (v1) {\lc100};
\node[rect] at (v2) {\lc110};
\node[rect] at (v3) {\lc011};
\node[rect] at (v4) {\lc111};
\node[rect] at (v5) {\lc101};
\node[rect] at (v6) {\lc111};
\draw[line width=1.2] (v2A)--(v2C);
\draw[line width=1.2] (v2B)--(v2D);
\draw[line width=1.2] (v5A)--(v5C);
\draw[line width=1.2] (v5B)--(v5D);
\end{tikzpicture}}
\quad+\quad
\scalebox{0.8}{\begin{tikzpicture}[vcenter]
\setupCoords{};
\draw[line width=1] (v2)--(v3)--(v4)--(v5)--(v6);
\node[rect] at (v1) {\lc001};
\node[rect] at (v2) {\lc110};
\node[rect] at (v3) {\lc011};
\node[rect] at (v4) {\lc111};
\node[rect] at (v5) {\lc101};
\node[rect] at (v6) {\lc111};
\draw[line width=1.2] (v4A)--(v4C);
\draw[line width=1.2] (v4B)--(v4D);
\end{tikzpicture}}
\end{equation}
To upper bound the two terms on the right-hand side of \eqref{eq:toy-split-v0}, we apply the induction hypothesis to the yield the following inequalities:
\begin{align}
\scalebox{0.8}{
\begin{tikzpicture}[vcenter]
\setupCoords{};
\draw[line width=1] (v2)--(v3)--(v4)--(v5)--(v6);
\node[rect] at (v1) {\lc001};
\node[rect] at (v2) {\lc110};
\node[rect] at (v3) {\lc011};
\node[rect] at (v4) {\lc111};
\node[rect] at (v5) {\lc101};
\node[rect] at (v6) {\lc110};
\end{tikzpicture}}
&\quad\le\quad 
\scalebox{0.8}{\begin{tikzpicture}[vcenter]
\moreCoords{};
\draw[line width=1] (v21)--(v25)--(v22);
\draw[line width=1] (v31)--(v34)--(v32)--(v33)--(v31);
\draw[line width=1] (v51)--(v55)--(v52);
\draw[line width=1] (v61)--(v64)--(v62)--(v63)--(v61);
\quad=\quad
\node[rect] at (v25) {\lc011};
\node[rect] at (v21) {\lc110};
\node[rect] at (v22) {\lc110};
\node[rect] at (v31) {\lc011};
\node[rect] at (v32) {\lc011};
\node[rect] at (v33) {\lc111};
\node[rect] at (v34) {\lc111};
\node[rect] at (v55) {\lc101};
\node[rect] at (v51) {\lc110};
\node[rect] at (v52) {\lc110};
\node[rect] at (v61) {\lc101};
\node[rect] at (v62) {\lc101};
\node[rect] at (v63) {\lc111};
\node[rect] at (v64) {\lc111};
\node at (p5) {$\frac{1}{2}$};
\node at (p3) {$\frac{1}{4}$};
\node at (p6) {$\frac{1}{4}$};
\node at (p2) {$\frac{1}{2}$};
\node at (p22) {$\color{white} \frac{1}{4}$};
\end{tikzpicture}}
\\
\scalebox{0.8}{
\begin{tikzpicture}[vcenter]
\setupCoords{};
\draw[line width=1] (v2)--(v3)--(v4)--(v5)--(v6);
\node[rect] at (v1) {\lc100};
\node[rect] at (v2) {\lc010};
\node[rect] at (v3) {\lc011};
\node[rect] at (v4) {\lc111};
\node[rect] at (v5) {\lc101};
\node[rect] at (v6) {\lc011};
\end{tikzpicture}}
&\quad\le\quad
\scalebox{0.8}{\begin{tikzpicture}[vcenter]
\moreCoords{};
\draw[line width=1] (v21)--(v25)--(v22);
\draw[line width=1] (v31)--(v34)--(v32)--(v33)--(v31);=
\draw[line width=1] (v51)--(v55)--(v52);
\draw[line width=1] (v61)--(v64)--(v62)--(v63)--(v61);
\node[rect] at (v25) {\lc011};
\node[rect] at (v21) {\lc010};
\node[rect] at (v22) {\lc010};
\node[rect] at (v31) {\lc011};
\node[rect] at (v32) {\lc011};
\node[rect] at (v33) {\lc111};
\node[rect] at (v34) {\lc111};
\node[rect] at (v55) {\lc101};
\node[rect] at (v51) {\lc011};
\node[rect] at (v52) {\lc011};
\node[rect] at (v61) {\lc101};
\node[rect] at (v62) {\lc101};
\node[rect] at (v63) {\lc111};
\node[rect] at (v64) {\lc111};
\node at (p5) {$\frac{1}{2}$};
\node at (p3) {$\frac{1}{4}$};
\node at (p6) {$\frac{1}{4}$};
\node at (p2) {$\frac{1}{2}$};
\node at (p22) {$\color{white} \frac{1}{4}$};
\end{tikzpicture}}
\end{align}

\textbf{First localization.} Here is the inequality that we are now left to prove:
\begin{equation}\label{eq:cancellation}
\scalebox{0.7}{
\begin{tikzpicture}[vcenter]
\moreCoords{};
\draw[line width=1] (v21)--(v25)--(v22);
\draw[line width=1,color=white!30!gray] (v31)--(v34)--(v32)--(v33)--(v31);
\draw[line width=1] (v51)--(v55)--(v52);
\draw[line width=1,color=white!30!gray] (v61)--(v64)--(v62)--(v63)--(v61);
\node[rect] at (v25) {\lc011};
\node[rect] at (v21) {\lc110};
\node[rect] at (v22) {\lc110};
\node[rect] at (v31) {\lcg011};
\node[rect] at (v32) {\lcg011};
\node[rect] at (v33) {\lcg111};
\node[rect] at (v34) {\lcg111};
\node[rect] at (v55) {\lc101};
\node[rect] at (v51) {\lc110};
\node[rect] at (v52) {\lc110};
\node[rect] at (v61) {\lcg101};
\node[rect] at (v62) {\lcg101};
\node[rect] at (v63) {\lcg111};
\node[rect] at (v64) {\lcg111};
\node at (p5) {$\frac{1}{2}$};
\node at (p3) {$\color{gray}\frac{1}{4}$};
\node at (p6) {$\color{gray}\frac{1}{4}$};
\node at (p2) {$\frac{1}{2}$};
\node at (p22) {$\color{white} \frac{1}{4}$};
\draw[line width=1] (x1)--(x2);
\draw[line width=1] (x3)--(x4);
\end{tikzpicture}}
\quad + \quad
\scalebox{0.7}{
\begin{tikzpicture}[vcenter]
\moreCoords{};
\draw[line width=1] (v21)--(v25)--(v22);
\draw[line width=1, color=white!30!gray] (v31)--(v34)--(v32)--(v33)--(v31);
\draw[line width=1] (v51)--(v55)--(v52);
\draw[line width=1,color=white!30!gray] (v61)--(v64)--(v62)--(v63)--(v61);
\node[rect] at (v25) {\lc011};
\node[rect] at (v21) {\lc010};
\node[rect] at (v22) {\lc010};
\node[rect] at (v31) {\lcg011};
\node[rect] at (v32) {\lcg011};
\node[rect] at (v33) {\lcg111};
\node[rect] at (v34) {\lcg111};
\node[rect] at (v55) {\lc101};
\node[rect] at (v51) {\lc011};
\node[rect] at (v52) {\lc011};
\node[rect] at (v61) {\lcg101};
\node[rect] at (v62) {\lcg101};
\node[rect] at (v63) {\lcg111};
\node[rect] at (v64) {\lcg111};
\node at (p5) {$\frac{1}{2}$};
\node at (p3) {$\color{gray}\frac{1}{4}$};
\node at (p6) {$\color{gray}\frac{1}{4}$};
\node at (p2) {$\frac{1}{2}$};
\node at (p22) {$\color{white} \frac{1}{4}$};
\draw[line width=1] (x1)--(x2);
\draw[line width=1] (x3)--(x4);
\end{tikzpicture}
}
\quad\le\quad
\scalebox{0.7}{
\begin{tikzpicture}[vcenter]
\moreCoords{};
\draw[line width=1] (v11)--(v14)--(v12)--(v13)--(v11);
\draw[line width=1] (v21)--(v24)--(v22)--(v23)--(v21);
\draw[line width=1,color=white!30!gray] (v31)--(v34)--(v32)--(v33)--(v31);
\draw[line width=1] (v41)--(v44)--(v42)--(v43)--(v41);
\draw[line width=1] (v51)--(v54)--(v52)--(v53)--(v51);
\draw[line width=1,color=white!30!gray] (v61)--(v64)--(v62)--(v63)--(v61);
\node[rect] at (v11) {\lc101};
\node[rect] at (v12) {\lc101};
\node[rect] at (v13) {\lc110};
\node[rect] at (v14) {\lc110};
\node[rect] at (v21) {\lc110};
\node[rect] at (v22) {\lc110};
\node[rect] at (v23) {\lc011};
\node[rect] at (v24) {\lc011};
\node[rect] at (v31) {\lcg011};
\node[rect] at (v32) {\lcg011};
\node[rect] at (v33) {\lcg111};
\node[rect] at (v34) {\lcg111};
\node[rect] at (v41) {\lc101};
\node[rect] at (v42) {\lc101};
\node[rect] at (v43) {\lc111};
\node[rect] at (v44) {\lc111};
\node[rect] at (v51) {\lc111};
\node[rect] at (v52) {\lc111};
\node[rect] at (v53) {\lc101};
\node[rect] at (v54) {\lc101};
\node[rect] at (v61) {\lcg101};
\node[rect] at (v62) {\lcg101};
\node[rect] at (v63) {\lcg111};
\node[rect] at (v64) {\lcg111};
\node at (p1) {$\frac{1}{4}$};
\node at (p2) {$\frac{1}{4}$};
\node at (p3) {$\color{gray}\frac{1}{4}$};
\node at (p4) {$\frac{1}{4}$};
\node at (p5) {$\frac{1}{4}$};
\node at (p6) {$\color{gray}\frac{1}{4}$};
\node at (p22) {$\color{white} \frac{1}{4}$};
\draw[line width=1] (x1)--(x2);
\draw[line width=1] (x3)--(x4);
\end{tikzpicture}
}
\end{equation}

Note that all the factors associated to edges that are more than 2 edges away from the deleted vertex in $G  = C_6$ are identical on both sides of the inequality, and thus they can be discarded (hence they are grayed and crossed out above).

We are left with showing the following inequality. Note that we have reduced the original inequality to a more local one involving only the edges of $G = C_6$ that are within two steps of the the deleted vertex (i.e., the edges of $G$ incident to a neighbor of the left-most vertex)
\begin{equation}
\begin{tikzpicture}[vcenter]
\moreCoords{};
\draw[line width=1] (v21B)--(v25B)--(v22B);
\draw[line width=1] (v51B)--(v55B)--(v52B);
\draw[line width=1] (c10)--(c1)--(c4)--(c40);
\draw[line width=1] (c20)--(c2)--(c3)--(c30);
\node[rect] at (v25B) {\lc011};
\node[rect] at (v21B) {\lc110};
\node[rect] at (v22B) {\lc110};
\node[rect] at (v55B) {\lc101};
\node[rect] at (v51B) {\lc110};
\node[rect] at (v52B) {\lc110};
\node at (p2B) {$\frac{1}{2}$};
\node at (p2BB) {$\color{white}\frac{1}{2}$};
\end{tikzpicture}
\quad +\quad 
\begin{tikzpicture}[vcenter]
\moreCoords{};
\draw[line width=1] (v21B)--(v25B)--(v22B);
\draw[line width=1] (v51B)--(v55B)--(v52B);
\draw[line width=1] (c10)--(c1)--(c4)--(c40);
\draw[line width=1] (c20)--(c2)--(c3)--(c30);
\node[rect] at (v25B) {\lc011};
\node[rect] at (v21B) {\lc010};
\node[rect] at (v22B) {\lc010};
\node[rect] at (v55B) {\lc101};
\node[rect] at (v51B) {\lc011};
\node[rect] at (v52B) {\lc011};
\node at (p2B) {$\frac{1}{2}$};
\node at (p2BB) {$\color{white}\frac{1}{2}$};
\end{tikzpicture}
\quad\le\quad
\begin{tikzpicture}[vcenter]
\draw[line width=1] (v11)--(v14)--(v12)--(v13)--(v11);
\draw[line width=1] (v21)--(v24)--(v22)--(v23)--(v21);
\draw[line width=1] (v41)--(v44)--(v42)--(v43)--(v41);
\draw[line width=1] (v51)--(v54)--(v52)--(v53)--(v51);
\node[rect] at (v11) {\lc101};
\node[rect] at (v12) {\lc101};
\node[rect] at (v13) {\lc110};
\node[rect] at (v14) {\lc110};
\node[rect] at (v21) {\lc110};
\node[rect] at (v22) {\lc110};
\node[rect] at (v23) {\lc011};
\node[rect] at (v24) {\lc011};
\node[rect] at (v41) {\lc101};
\node[rect] at (v42) {\lc101};
\node[rect] at (v43) {\lc111};
\node[rect] at (v44) {\lc111};
\node[rect] at (v51) {\lc111};
\node[rect] at (v52) {\lc111};
\node[rect] at (v53) {\lc101};
\node[rect] at (v54) {\lc101};
\node at (p1) {$\frac{1}{4}$};
\node at (p2) {$\frac{1}{4}$};
\node at (p4) {$\frac{1}{4}$};
\node at (p5) {$\frac{1}{4}$};
\node at (p22) {$\color{white} \frac{1}{4}$};
\end{tikzpicture}
\end{equation}

\textbf{Second localization.} Now let us apply the Cauchy--Schwarz inequality in the form of $\sqrt{a_1 b_1} + \sqrt{a_2b_2} \le \sqrt{a_1+a_2}\sqrt{b_1+b_2}$ to the left-hand side above. The remaining inequality to show follows by taking the product of the following two inequalities (corresponding to the top and bottom halves of the above inequality after Cauchy--Schwarz):
\begin{align}
\label{eq:toy-post-cs-top}
\begin{tikzpicture}[vcenter]
\draw[line width=1] (v21B)--(v25B)--(v22B);
\node[rect] at (v25B) {\lc011};
\node[rect] at (v21B) {\lc110};
\node[rect] at (v22B) {\lc110};
\end{tikzpicture}
&\quad+\quad 
\begin{tikzpicture}[vcenter]
\draw[line width=1] (v21B)--(v25B)--(v22B);
\node[rect] at (v25B) {\lc011};
\node[rect] at (v21B) {\lc010};
\node[rect] at (v22B) {\lc010};
\end{tikzpicture}
\quad\le \quad 
\begin{tikzpicture}[vcenter]
\draw[line width=1] (v11)--(v14)--(v12)--(v13)--(v11);
\draw[line width=1] (v21C)--(v24C)--(v22C)--(v23C)--(v21C);
\node[rect] at (v11) {\lc101};
\node[rect] at (v12) {\lc101};
\node[rect] at (v13) {\lc110};
\node[rect] at (v14) {\lc110};
\node[rect] at (v21C) {\lc110};
\node[rect] at (v22C) {\lc110};
\node[rect] at (v23C) {\lc011};
\node[rect] at (v24C) {\lc011};
\node at (p1) {$\frac{1}{2}$};
\node at (p2C) {$\frac{1}{2}$};
\node at (p2CCC) {$\color{white}\frac{1}{2}$};
\end{tikzpicture}
\\
\label{eq:toy-post-cs-bottom}
\begin{tikzpicture}[vcenter]
\draw[line width=1] (v51B)--(v55B)--(v52B);
\node[rect] at (v55B) {\lc101};
\node[rect] at (v51B) {\lc110};
\node[rect] at (v52B) {\lc110};
\end{tikzpicture}
&\quad+\quad
\begin{tikzpicture}[vcenter]
\draw[line width=1] (v51B)--(v55B)--(v52B);
\node[rect] at (v55B) {\lc101};
\node[rect] at (v51B) {\lc011};
\node[rect] at (v52B) {\lc011};
\end{tikzpicture}
\quad\le \quad 
\begin{tikzpicture}[vcenter]
\draw[line width=1] (v41)--(v44)--(v42)--(v43)--(v41);
\draw[line width=1] (v51C)--(v54C)--(v52C)--(v53C)--(v51C);
\node[rect] at (v41) {\lc101};
\node[rect] at (v42) {\lc101};
\node[rect] at (v43) {\lc111};
\node[rect] at (v44) {\lc111};
\node[rect] at (v51C) {\lc111};
\node[rect] at (v52C) {\lc111};
\node[rect] at (v53C) {\lc101};
\node[rect] at (v54C) {\lc101};
\node at (p4) {$\frac{1}{2}$};
\node at (p5C) {$\frac{1}{2}$};
\node at (p5CC) {$\color{white}\frac{1}{2}$};
\end{tikzpicture}
\end{align}

Note that each inequality now involves only a two-edge path in $G = C_6$ starting from the deleted vertex. Thus we have further localized the inequality that we wish to prove.

Let us explain how to prove \eqref{eq:toy-post-cs-top}, as the proof of \eqref{eq:toy-post-cs-bottom} is analogous (in this specific example \eqref{eq:toy-post-cs-bottom} is actually an equality). As in \eqref{eq:toy-split-v0}, the left-hand side of \eqref{eq:toy-post-cs-top} can be rewritten as:
\begin{equation}\label{eq:toy-c4-equal}
\begin{tikzpicture}[vcenter]
\draw[line width=1] (v51D)--(v54D)--(v52D)--(v53D)--(v51D);

\node[rect] at (v51D) {\lc110};
\node[rect] at (v52D) {\lc110};
\node[rect] at (v53D) {\lc101};
\node[rect] at (v54D) {\lc011};

\end{tikzpicture}
\quad=\quad
\begin{tikzpicture}[vcenter]
\draw[line width=1] (v21B)--(v25B)--(v22B);
\node[rect] at (v25B) {\lc011};
\node[rect] at (v21B) {\lc110};
\node[rect] at (v22B) {\lc110};
\end{tikzpicture}
\quad+\quad
\begin{tikzpicture}[vcenter]
\draw[line width=1] (v21B)--(v25B)--(v22B);
\node[rect] at (v25B) {\lc011};
\node[rect] at (v21B) {\lc010};
\node[rect] at (v22B) {\lc010};
\end{tikzpicture}
\end{equation}

Now we are left with proving the following inequality. Note that it has the same form as the strong induction hypothesis (e.g., \eqref{eq:toy-c6}). While it may be tempting to quote the induction hypothesis, it is instead more faithful to view this as a version of the local inequality which we need to prove. 
\begin{equation} \label{eq:toy-final-C4}
\begin{tikzpicture}[vcenter]
\draw[line width=1] (v51D)--(v54D)--(v52D)--(v53D)--(v51D);

\node[rect] at (v51D) {\lc110};
\node[rect] at (v52D) {\lc110};
\node[rect] at (v53D) {\lc101};
\node[rect] at (v54D) {\lc011};

\end{tikzpicture}
\quad\le \quad
\begin{tikzpicture}[vcenter]
\draw[line width=1] (v51C)--(v54C)--(v52C)--(v53C)--(v51C);

\node[rect] at (v53C) {\lc110};
\node[rect] at (v54C) {\lc110};
\node[rect] at (v51C) {\lc101};
\node[rect] at (v52C) {\lc101};
\node at (p5C) {$\frac{1}{2}$};
\node at (p5CC) {$\color{white} \frac{1}{2}$};
\end{tikzpicture}
\quad 
\begin{tikzpicture}[vcenter]
\draw[line width=1] (v51C)--(v54C)--(v52C)--(v53C)--(v51C);
\node[rect] at (v51C) {\lc110};
\node[rect] at (v52C) {\lc110};
\node[rect] at (v53C) {\lc011};
\node[rect] at (v54C) {\lc011};
\node at (p5C) {$\frac{1}{2}$};
\node at (p5CC) {$\color{white} \frac{1}{2}$};
\end{tikzpicture}
\end{equation}
This inequality follows from the Cauchy-Schwarz inequality applied as follows:
\begin{align*}
&\quad\sum_{a,b_1,b_2,c}  f(a,b_1)f(a, b_2)g(b_1,c)g(b_2,c)
\\ &=\sum_{b_1,b_2} \left(\sum_af(a,b_1)f(a, b_2)\right)\left(\sum_cg(b_1,c)g(b_2,c)\right)
\\ &\le \left(\sum_{b_1,b_2}\left(\sum_{a}f(a,b_1)f(a, b_2)\right)^2 \right)^{\frac{1}{2}}\left(\sum_{b_1,b_2}\left(\sum_{c}g(b_1,c)g(b_2,c)\right)^2 \right)^{\frac{1}{2}}
\\ &= \left(\sum_{a_1,a_2,b_1,b_2}f(a_1,b_1)f(a_1,b_2)f(a_2,b_1)f(a_2,b_2)\right)^{\frac{1}{2}}\left(\sum_{b_1,b_2,c_1,c_2}g(b_1,c_1)g(b_2,c_1)g(b_1,c_2)g(b_2,c_2)\right)^{\frac{1}{2}}.
\end{align*}
Here the variables $a, b_1, \dots$ range over $\{\text{red}, \text{green}, \text{blue}\}$ in the sums, $f(a,b)$ is the indicator function associated to coloring the top-left edge on the left-hand side of \eqref{eq:toy-final-C4}, i.e.,
\[
f(a,b) = \begin{cases}
1 & \text{if } a \in \{\text{red}, \text{blue}\}, b \in \{\text{red}, \text{green}\}, \text{ and } a \ne b, \\
0 & \text{otherwise,}
\end{cases}
\]
and $g(b,c)$ is analogously defined for the top-right edge on the left-hand side of \eqref{eq:toy-final-C4}. This completes the proof of \eqref{eq:toy-post-cs-top}.

\subsection*{Further complications}
In the general setting, the first and second localization steps are analogous to the toy calculation above. In particular, the induction proceeds by first selecting the color of a maximum degree vertex $w$ (in the above calculation, the left-most vertex), and then updating the lists of colors in the neighborhood of $w$ for each color selection. We then apply the inductive hypothesis followed by the first localization as in \eqref{eq:cancellation}, reducing the problem to just considering edges in a radius 2 neighborhood of $w$. The second localization is in general an application of H\"older's inequality, which reduces the problem to inequalities on two-edge paths and triangles. The analysis is somewhat easier in the triangle-free case, which is done in \cref{sec:triangle-free}, as we only need to prove one type of local inequality. When $G$ contains triangles, which is done in \cref{sec:color} for colorings, the presence of triangles require the analysis of additional local inequalities that are more difficult to handle. The additional difficulty is expected since the local inequalities involving triangles cannot be true for all targets/models $H$, so the proofs need to use more specific knowledge of the model.

Even in the triangle-free case, the local inequalities for two-edge paths are in general more involved than shown above in \eqref{eq:toy-post-cs-top} and \eqref{eq:toy-post-cs-bottom}. This is because the equality \eqref{eq:toy-c4-equal} turns out to be coincidental to the graph $G=C_6$; in general, one side of these inequalities is a summation of terms with fractional exponents. We handle this difficulty by defining an interpolation between the terms of the local inequality, and proving log-convexity with respect to the underlying parameter of this interpolation.

In \cref{sec:color}, we handle local inequalities for triangles in the case where $H$ is the complete graph, possibly with loops on some of its vertices. In this particular case, the result follows from additional log-convexity results along with an intricate analysis of correlation inequalities on symmetric polynomials. 

For the clique-maximization result \cref{thm:clique-max}, proved in \cref{sec:clique}, the spirit of the solution is similar, although the execution differs. In this case, since the upper bound is a product over vertices, the analogous localization steps result in statements for each \emph{vertex} in the neighborhood of $w$. Here again we have fractional exponents in general, which are handled with a series of several interpolations, each of which is shown to be log-convex. 

\section{Inequality on triangle-free graphs} \label{sec:triangle-free}

The goal of this section is to prove \cref{thm:graph-BL}. 

\subsection{Some preliminary inequalities}

In the lemmas, we omit stating the obvious integrability hypotheses.

\begin{lemma} \label{lem:mixed-norm}
	For nonnegative functions $g(s,u)$ and $h(s,v)$, and real $q \ge 1$, one has
	\begin{multline*}
		\left(\int \left(\int  g(s,u)h(s,v)\, du ds\right)^q dv\right)^{2/q} \\
	\le \left( \int g(s,u)g(s,u') \, dudu'ds\right) \left( \int \left( \int  h(s,v)h(s,v')  \, ds \right )^q dvdv' \right)^{1/q}.
	\end{multline*}
\end{lemma}

\begin{proof}
	Let $1/q + 1/q' = 1$. We have
	\begin{align*}
	\text{LHS} &= \sup_{\norm{f}_{L^{q'}} \le 1} \left( \int g(s,u) h(s,v) f(v) \,dudvds\right)^2 && \text{\small [$L^p$ duality]} \\
		&\le \sup_{\norm{f}_{L^{q'}} \le 1} \left( \int \left(\int g(s,u) \, du\right)^2  \,ds \right)  \left( \int \left( \int h(s,v) f(v)  dv\right)^2 \,ds\right) && \text{\small [Cauchy--Schwarz]} \\
		&= \left( \int g(s,u)g(s,u') \, dudu'ds\right) \sup_{\norm{f}_{L^{q'}} \le 1}  \int h(s,v)h(s,v')f(v)f(v') \,dvdv' \,ds,
	\end{align*}
	which is at most the RHS, by using $L^p$ duality once again, noting that the function $(v,v') \mapsto f(v)f(v')$ has $L^{q'}$ norm $\norm{f}_{L^{q'}}^2 \le 1$.
\end{proof}

\begin{remark} Define the mixed $L^{p,q}$ matrix norm of $A = (a_{ij})$ by 
\[
\norm{A}_{L_{p,q}} := \left(\sum_i \left( \sum_j \abs{a_{ij}}^p \right)^{q/p} \right)^{1/q}.
\]
\cref{lem:mixed-norm} is equivalent to the following inequality: If $q \ge 1$, $A \in \RR_{\ge 0}^{m\times n}$ and $B \in \RR_{\ge 0}^{m \times k}$, then 
\begin{equation} \label{eq:matrix-mixed-norm-ineq}
\norm{A^\T B}_{L_{1, q}}^2 \le \norm{A^\T A}_{L_{q, q}}\norm{B^\T B}_{L_{1, 1}}.
\end{equation}
Here is the above proof written out in the language of matrices and vectors, with $\bm 1_m$ denoting the all ones vector:
\begin{multline*}
\norm{A^\T B}^2_{L_{1, q}}
= \norm{A^\T B\bm 1_m}_q^2
= \sup_{\norm{\bm u}_{q'}\le 1}\langle \bm u, A^\T B\bm 1_m\rangle^2
= \sup_{\norm{\bm u}_{q'}\le 1}\langle A \bm u, B \bm 1_m\rangle^2
\\ \le \sup_{\norm{\bm u}_{q'}\le 1}\langle A \bm u, A \bm u\rangle \langle B \bm 1_m, B\bm 1_m\rangle
= \sup_{\norm{\bm u}_{q'}\le 1} \langle \bm u\bm u^\T, A^\T A\rangle \langle \bm 1_m \bm 1_m^\T, B^\T B\rangle 
\le \norm{A^\T A}_{L_{q, q}} \norm{B^\T B}_{L_{1, 1}}.
\end{multline*}
We do not know if the inequality can be extended to $\norm{A^\T B}_{L_{p, q}}^2 \le \norm{A^\T A}_{L_{q, q}}\norm{B^\T B}_{L_{p, p}}$ for all reals $1 \le p \le q$ (this is true for positive integer $p$ by a tensor-power argument). Note that the above proof holds for all real matrices $A$ and $B$ provided that $A^\T A, B^\T B, A^\T B$ have nonnegative entries.
\end{remark}

\begin{lemma} \label{lem:mixed-norm-2}
For nonnegative functions $f(s,t)$, $g(s,t,u)$, $h(s,t,v)$, and real $q \ge 1$, 
\begin{multline*}
\left( \int \left(\int  f(s,t)g(s,t,u)h(s,t,v)\, du ds\right)^q  dt dv \right)^2
\\ 
\le
\left( \int \left( \int f(s,t)g(s,t,u)g(s,t,u') \, dudu'ds\right)^q \,dt \right)
\left( \int \left( \int f(s,t) h(s,t,v)h(s,t,v')  \, ds \right )^q dt dvdv' \right) 
\end{multline*}
\end{lemma}

\begin{proof}
	By replacing $g(s,t,u)$ by $f(s,t)^{1/2}g(s,t,u)$ and $h(s,t,v)$ by $f(s,t)^{1/2}h(s,t,v)$, we may assume that $f = 1$. By the Cauchy--Schwarz inequality with respect to $dt$, the right-hand side is at least
	\[
	\left( \int \left( \int g(s,t,u)g(s,t,u') \, dudu'ds\right)^{q/2} \left( \int \left( \int  h(s,t,v)h(s,t,v')  \, ds \right )^q dvdv' \right)^{1/2} dt  \right)^2.
	\]
	It suffices to show that, for every fixed $t$, one has
	\begin{multline*}	
	\int \left(\int  g(s,t,u)h(s,t,v)\, du ds\right)^q  dv \\
	\le \left( \int g(s,t,u)g(s,t,u') \, dudu'ds\right)^{q/2} \left( \int \left( \int  h(s,t,v)h(s,t,v')  \, ds \right )^q dvdv' \right)^{1/2},
	\end{multline*}
	which is \cref{lem:mixed-norm} applied to the functions $g(s,u) = g(s,t,u)$ and $h(s,v) = h(s,t,v)$.
\end{proof}

The following lemma is a ``local'' inequality that the proof of \cref{thm:graph-BL} will reduce to.

\begin{lemma}[Local inequality] \label{lem:123}
	Let $f_{12} \colon \Omega_1 \times \Omega_2 \to \RR_{\ge 0}$ and $f_{23} \colon \Omega_2 \times \Omega_3 \to \RR_{\ge 0}$ be measurable functions. Let $1\le \beta\le \Delta$ and $\gamma \ge 2$ be integers. For $x \in \Omega_1$, define $f_{23}^x \colon \Omega_2 \times \Omega_3 \to \RR_{\ge 0}$ by $f_{23}^x(y,z) := f_{12}(x,y)^{1/(\gamma-1)}f_{23}(y,z)$ and $\norm{\cdot}_{K_{a,b}}$ as in \cref{defn:graph-norm}. Then
	\[
		\int_{\Omega_1} \norm{f_{23}^x}_{K_{\beta,\gamma-1}}^{\Delta(\gamma-1)} \, dx\le  \norm{f_{12}}_{K_{\gamma,\Delta}}^\Delta \norm{f_{23}}_{K_{\beta, \gamma}}^{\Delta(\gamma-1)}.
	\]
\end{lemma}

\begin{proof}
Define, for nonnegative integers $a, b, c$,
\[
M_{a,b,c} = 
\int_{\Omega_1^a}
	\left( 
		\int_{\Omega_2^b \times \Omega_3^c}  
			\prod_{i \in [a], j\in[b]} f_{12} (x^{(i)}, y^{(j)}) 
			\prod_{j\in[b], k \in [c]} f_{23} (y^{(j)}, z^{(k)}) 
		\, d\bm y^{[b]} d \bm z^{[c]}
	\right)^{\Delta/\beta}
\, d \bm x^{[a]},
\]
where $d\bm x^{[a]} = d x^{(1)} \dotsm d x^{(a)}$, $d\bm y^{[b]} = d y^{(1)} \dotsm d y^{(b)}$, and $d\bm z^{[c]} = d z^{(1)} \dotsm d z^{(c)}$.
By expanding $f_{23}^{x}$, we have
\[
\int_{\Omega_1} \norm{f_{23}^{x}}_{K_{\beta,\gamma-1}}^{\Delta(\gamma-1)} \, d x  = M_{1,\beta,\gamma-1}.
\]
Also
\[
\|f_{12}\|_{K_{\gamma,\Delta}}^{\gamma\Delta} = M_{\gamma, \beta, 0}
\quad \text{and} \quad
\|f_{23}\|_{K_{\beta,\gamma}}^{\beta \gamma} = M_{0, \beta,\gamma}^{\beta/\Delta}.
\]
Thus the claimed inequality can be written as
\[
M_{1,\beta,\gamma-1} \le  M_{\gamma, \beta, 0}^{1/\gamma}  M_{0, \beta,\gamma}^{1- 1/\gamma},
\]
which would follow from $(M_{i,\beta,\gamma-i})_{0\le i\le \gamma}$ being log-convex. Thus it suffices to prove that 
\[
M_{a+1,b,c+1}^2 \le M_{a+2,b,c} M_{a,b,c+2},
\]
for all nonnegative integers $a,b,c$. This inequality follows from \cref{lem:mixed-norm-2}, after setting $q = \Delta/\beta\ge 1$, with
\begin{align*}
	\bm s &= (y^{(1)}, \dotsc, y^{(b)}, z^{(1)}, \dotsc, z^{(c)}), \\
	\bm t &= (x^{(1)}, \dotsc, x^{(a)}), \\
	\bm u &= (z^{(c+1)}), \quad \bm u' = (z^{(c+2)}), \\	
	\bm v &= (x^{(a+1)}), \quad \bm v' = (x^{(a+2)}), \\
	f(\bm s, \bm t)
	&= \prod_{i \in [a], j\in[b]} f_{12} (x^{(i)}, y^{(j)}) 
			\prod_{j\in[b], k \in [c]} f_{23} (y^{(j)}, z^{(k)}) ,
	\\
	g(\bm s, \bm t, \bm u)
     &= \prod_{j \in [b]} f_{23} (y^{(j)}, z^{(c+1)} ), \qquad 
	g(\bm s, \bm t, \bm u')
	= \prod_{j \in [b]} f_{23} (y^{(j)}, z^{(c+2)} )
	\\
	h(\bm s, \bm t, \bm v)
	&= \prod_{j \in [b]} f_{12} (x^{(a+1)} ,y^{(j)}), \qquad
	 h(\bm s, \bm t, \bm v')
	= \prod_{j \in [b]} f_{12} (x^{(a+2)} ,y^{(j)}). \qedhere
\end{align*}
\end{proof}

\subsection{Proof of Theorem~\ref{thm:graph-BL}} \label{sec:pf-BL-main}

We apply induction on the number of vertices in $G = (V,E)$. The base case $\abs{E} = 0$ is trivial. Let $\Delta$ be the maximum degree of $G$, and let $w$ be a vertex of degree $\Delta$ in $G$. The idea of the following calculation is to consider what happens when we condition on a certain color (i.e., element of $\Omega_w$) assigned to $w$.

\smallskip

\noindent\emph{Notation.} For $k \in \{0, 1, 2, \dots\}\cup\{\infty\}$, let $V_k$ be the set of vertices at distance exactly $k$ from $w$. For $0 \le i \le j \le i+1$, let $E_{ij}$ be the edges with one endpoint in $V_i$ and the other in $V_j$. Let $V_{\ge k} = \bigcup_{i \ge k} V_i$, $E_{\ge k} = \bigcup_{k \le i \le j} E_{ij}$, and $E_{> k} = E_{k,k+1} \cup E_{\ge k+1}$. Note that $V = V_{\ge 0}$ and $E = E_{> 0}$.

Let $I_1$ be the set of vertices in $V_1$ whose neighborhood is exactly $\{w\}$.

Although we treat edges as unordered pairs, when we write $vu \in E_{ij}$, we always mean $v \in V_i$ and $u \in V_j$. On the other hand, when we range over $uv \in E_{ii}$, we do not count $uv$ and $vu$ separately.

\begin{figure}[t]\centering
\includegraphics{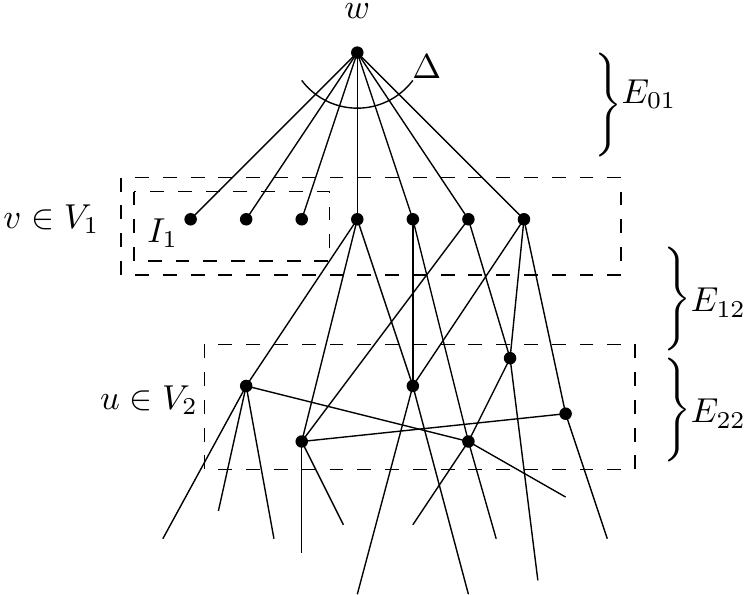}	
\caption{Labels of vertices and edges in the proof of \cref{thm:graph-BL}.}
\end{figure}

For any $S \subseteq V$, write $\Omega_S := \prod_{v \in S} \Omega_v$ and $d \bm x_S := \prod_{v \in S} dx_v$.

For $vu \in E_{12}$ with $v \in V_1$ and $u \in V_2$, and $x_w \in \Omega_w$, define $f_{vu}^{x_w} \colon \Omega_v \times \Omega_u \to \RR$ by
\[
f_{vu}^{x_w}(x_v,x_u) := f_{wv}(x_w,x_v)^{1/(d_v - 1)} f_{vu}(x_v, x_u).
\]

\smallskip

By distributing the $E_{01}$ factors to $E_{12}$, we have
\begin{align*}
	&\int_{\Omega_V} \prod_{vu\in E} f_{vu}(x_v,x_u) \, d\bm x_V
	\\ 
	&= \int_{\Omega_{V\setminus I_1}} \prod_{v \in I_1} \left(\int_{\Omega_v} f_{wv}(x_w,x_v) \, dx_v\right) \prod_{vu\in E_{12}} \left( f_{wv}(x_w,x_v)^{1/(d_v - 1)} f_{vu}(x_v,x_u) \right) \prod_{vu\in E_{\ge 2}} f_{vu}(x_v,x_u) \,
	 d\bm x_{V\setminus I_1}
	\\
	&\le \int_{\Omega_w} \prod_{v \in I_1} \norm{f_{wv}(x_w,\cdot)}_1 \prod_{vu\in E_{12}} \norm{f_{vu}^{x_w}}_{K_{d_u,d_v-1}} \prod_{vu\in E_{\ge 2}} \norm{f_{vu}}_{K_{d_u,d_v}} \, d x_w,
\end{align*}
where in the last step we applied the induction hypothesis to $G-w$ (the graph $G$ with the vertex $w$ removed along with all its incident edges).

Now fix $x_w \in \Omega_w$. It remains to prove the bound
\[
\int_{\Omega_w} \prod_{v \in I_1} \norm{f_{wv}(x_w,\cdot)}_1 \prod_{vu\in E_{12}} \norm{f_{vu}^{x_w}}_{K_{d_u,d_v-1}} \prod_{vu\in E_{\ge 2}} \norm{f_{vu}}_{K_{d_u,d_v}} \, d x_w
\le \prod_{vu\in E} \|f_{vu}\|_{K_{d_u,d_v}}.
 \]
 \textbf{First localization.}
Observing that the factor $\norm{f_{vu}}_{K_{d_u,d_v}}$ appears on both sides whenever $vu \in E_{\ge 2}$, we see that it suffices to prove
\[
\int_{\Omega_w} \prod_{v \in I_1} \norm{f_{wv}(x_w,\cdot)}_1 \prod_{vu\in E_{12}} \norm{f_{vu}^{x_w}}_{K_{d_u,d_v-1}} \, d x_w
\le \prod_{vu\in E_{01} \cup E_{12}} \|f_{vu}\|_{K_{d_u,d_v}}.
\]
By distributing the $E_{01}$ factors on the RHS to $E_{12}$, we can rewrite the above inequality as
\begin{equation}
	\label{eq:after-trim}
\int_{\Omega_w} \prod_{v \in I_1} \norm{f_{wv}(x_w,\cdot)}_1 \prod_{vu\in E_{12}} \norm{f_{vu}^{x_w}}_{K_{d_u,d_v-1}} \, d x_w
\le \prod_{v \in I_1} \|f_{wv}\|_{K_{1,\Delta}} \prod_{vu \in E_{12}} \left( \|f_{wv}\|_{K_{d_v,\Delta}}^{1/(d_v-1)} \|f_{vu}\|_{K_{d_u,d_v}}\right).
\end{equation}
From now on until the end of the proof, by convention, we use the letter $v$ to denote a vertex in $V_1$ and $u$ for a vertex in $V_2$.

\noindent \textbf{Second localization.}
Applying H\"older's inequality with exponents given by the summands of 
\[
\sum_{v \in I_1} \frac{1}{\Delta} + \sum_{vu \in E_{12}} \frac{1}{\Delta(d_v-1)} = 1,
\]
we upper bound the left-hand side of \cref{eq:after-trim} by
\begin{multline*}
\int_{\Omega_w} \prod_{v \in I_1} \norm{f_{wv}(x_w,\cdot)}_1 \prod_{uv\in E_{12}} \norm{f_{vu}^{x_w}}_{K_{d_u,d_v-1}} \, d x_w\
\\
\le 
\prod_{v \in I_1} \left( \int_{\Omega_w} \norm{f_{wv}(x_w,\cdot)}_1^\Delta \, d x_w\right)^{\frac{1}{\Delta}}
\prod_{vu\in E_{12}} \left( \int_{\Omega_w} \norm{f_{vu}^{x_w}}_{K_{d_u,d_v-1}}^{\Delta(d_v-1)} \, d x_w\right)^{\frac{1}
{\Delta(d_v-1)}}.
\end{multline*}
Comparing with the desired right-hand side of \eqref{eq:after-trim}, we have
\[
\int_{\Omega_w} \norm{f_{wv}(x_w,\cdot)}_1^\Delta \, d x_w = \norm{f_{wv}}_{K_{1,\Delta}}^\Delta
\]
and, by \cref{lem:123},  the local inequality, for every $vu \in E_{12}$,
\[
	\int_{\Omega_w} \norm{f_{vu}^{x_w}}_{K_{d_u,d_v-1}}^{\Delta(d_v-1)} \, d x_w  
	\le \|f_{wv}\|_{K_{d_v,\Delta}}^\Delta  \|f_{vu}\|_{K_{d_u,d_v}}^{\Delta(d_v-1)},
\]
which proves \eqref{eq:after-trim}. This concludes the proof of \cref{thm:graph-BL}. \qed

\subsection{Necessity of the triangle-free hypothesis}

Now we prove \cref{prop:counterexample-triangle}, showing that triangle-free hypotheses on $G$ in \cref{thm:K3-free-hom-irreg} (and hence also \cref{thm:graph-BL}) cannot be removed. 

\begin{proof}[Proof of \cref{prop:counterexample-triangle}]
	Let $H_\epsilon$ be a weighted graph on two vertices each with vertex weight $1/2$, and edge-weight ``adjacency'' matrix
	\[
	\begin{pmatrix}
		1+2\epsilon & 1 \\
		1 & 1+2\epsilon
	\end{pmatrix},
	\]
	i.e., a loop with weight $1+2\epsilon$ on each vertex, and an edge of unit weight between the two vertices. For every graph $G$, one has, for small $\epsilon$,
	\begin{align*}
	\hom(G, H_\epsilon)
	&= \EE_{\bm x \in \{0,1\}^{V(G)}} (1 + 2\epsilon)^{\abs{\set{uv \in E(G): x_u = x_v}}} 
	\\
	&= 1 + \abs{E(G)} \epsilon + \binom{\abs{E(G)}}{2} \epsilon^2 + \left(\binom{\abs{E(G)}}{3} + \abs{T(G)}\right) \epsilon^3 + O(\epsilon^4)
	\\
	&= \paren{ 1+ \epsilon}^{\abs{E(G)}} + \abs{T(G)}\epsilon^3 + O(\epsilon^4),
	\end{align*}
	where $T(G)$ is the set of triangles in $G$. Indeed, the coefficient of $\epsilon^k$ comes from examining each $k$-edge subset of $E(G)$ and determining the probability that each connected component of this $k$-edge subset receives the same color in $\bm x$. Thus,
	\[
	\hom(G, H_\epsilon)^{1/\abs{E(G)}} = 1 + \epsilon + \frac{\abs{T(G)}}{\abs{E(G)}} \epsilon^3 + O(\epsilon^4).
	\]
	On the other hand, since $K_{a,b}$ is always triangle-free,
	\[
	\prod_{uv \in E(G)} \hom(K_{d_u,d_v}, H)^{1/(d_ud_v)} 
	= (1 + \epsilon + O(\epsilon^4))^{\abs{E(G)}}.
	\]
	Thus, for every graph $G$ with at least one triangle, for sufficiently small $\epsilon > 0$, the weighted graph $H = H_\epsilon$ satisfies
    \begin{equation}
		\label{eq:triangle-counter}
			\hom(G, H) > \prod_{uv \in E(G)} \hom(K_{d_u,d_v}, H)^{1/(d_ud_v)}.
	\end{equation}
	Finally, one can obtain a simple graph $H$ satisfying \cref{eq:triangle-counter} by a standard graph limit argument (e.g., \cite{BCLS08} \cite[Ch.~10]{Lov}). Indeed, note that \cref{eq:triangle-counter} is also true for $H = H_\epsilon/2$ (obtained by scaling all edge weights of $H_\epsilon$ by $1/2$), since scaling the edge weights of $H$ by a factor of 2 changes both sides of \cref{eq:triangle-counter} by the same factor $2^{|E(G)|}$. Now take a sequence of simple graphs $H_n$ converging to $H_\epsilon/2$ in the sense that $\hom(F, H_n)/|V(H_n)|^{|V(F)|} \to \hom(F, H_\epsilon/2)$ as $n \to \infty$ for every graph $F$. Then for sufficiently large $n$, the simple graph $H = H_n$ satisfies \cref{eq:triangle-counter}, as both sides are scaled by $|V(H)|$ raised to the same exponent $|V(G)| = \sum_{uv \in E(G)} (d_u + d_v)/(d_u d_v)$.
	\end{proof}

\subsection{Antiferromagnetic 2-spin models are biclique-maximizing} \label{sec:2spin}

Here we prove the part of the claim in \cref{cor:2spin} that every 2-spin antiferromagnetic model is biclique-maximizing. The proof relies on the bipartite swapping trick introduced in \cite{Zhao10} for \cref{q:ind} on upper bounding the number of independent sets. It was proved that $i(G)^2 \le i(G \times K_2)$, thereby reducing \cref{q:ind} to the bipartite setting, which had been resolved earlier~\cite{Kahn01}. The following inequality is a weighted generalization of the same inequality for independent sets. The proof is adapted from \cite{Zhao10}, and we include it here for completeness. It can also be extended further as in \cite{Zhao11} to a larger class of weighted $H$, though we omit the details. See \cite[Theorem~1.13]{Csi16ar} for another proof based on similar ideas. 

\begin{lemma} \label{lem:bst}
	Let $G$ be a graph and $H$ be a 2-spin antiferromagnetic model. Then
	\begin{equation} \label{eq:bst}
		\hom(G, H)^2 \le \hom(G \times K_2, H).
	\end{equation}
\end{lemma}

Since $G \times K_2$ is bipartite, \cref{lem:bst} followed by \cref{thm:graph-BL} (or \cref{thm:K3-free-hom-irreg} for weighted $H$) gives
\[
\hom(G, H) \le \hom(G \times K_2, H)^{1/2} \le  \prod_{uv \in E(G)} \hom(K_{d_u,d_v}, H)^{1/(d_ud_v)},
\]
as the degree-degree distribution does not change when $G$ is lifted to $G \times K_2$. 

It remains open whether \eqref{eq:bst} holds for $H = K_q$, corresponding to proper colorings.

\begin{proof}[Proof of \cref{lem:bst}]
	Let $G = (V,E)$. Let $\Omega = \{0,1\}$ be a two-point measure space and let $H \colon \Omega \times \Omega \to \RR_{\ge 0}$ be antiferromagnetic, or equivalently, $H(0,0)H(1,1) \le H(0,1)^2$. Let 
	\[
	  S_= = \{(\bm x, \bm y, \bm z) :  z_{uv} \leq H(x_u, x_v) H(y_u, y_v) \ \forall uv \in E\} \subseteq \Omega^{V} \times \Omega^{V} \times \RR_{\ge 0}^{E}
	\]
	and
	\[
	  S_\times = \{(\bm x, \bm y, \bm z) :  z_{uv} \leq H(x_u, y_v) H(y_u, x_v) \ \forall uv \in E\}
	   \subseteq \Omega^{V} \times \Omega^{V} \times \RR_{\ge 0}^{E},
	\]
	where $\bm x = (x_{v})_{v \in V(G)} \in \Omega^{V}$, $\bm y = (y_v)_{v \in V} \in \Omega^{V}$ and $\bm z = (z_e)_{e \in E} \in \RR_{\ge 0}^{E}$.
	Note that $\hom(G, H)^2$ is equal to the measure of $S_=$, and $\hom(G \times K_2, H)$ equals to the measure of $S_\times$. Thus the lemma reduces to constructing a measure-preserving injection $\phi \colon S_= \to S_\times$.
	
	For any $(\bm x, \bm y, \bm z) \in \Omega^V \times \Omega^V \times \RR_{\ge 0}^E$, say that an edge $uv \in E$ is \emph{unsafe} with respect to $(\bm x, \bm y, \bm z)$ if either $z_{uv} > H(x_u, x_v) H(y_u,y_v)$ or $z_{uv} > H(x_u,y_v) H(y_u,x_v)$.  Fixing $(\bm x, \bm y, \bm z) \in S_=$, if $uv$ is unsafe, then $H(x_u, y_v) H(y_u, x_v) < z_{uv} \le H(x_u, x_v) H(y_u, y_v)$ (the former due to being unsafe, and the latter due to the definition of $S_=$). Recall $\Omega = \{0,1\}$. Since $H$ is 2-spin antiferromagnetic, the only way to satisfy $H(x_u, y_v) H(y_u, x_v) < H(x_u, x_v) H(y_u, y_v)$ is that one of the endpoints of $uv$, say $u$, has $(x_u, y_u) = (0,1)$, and the other endpoint $v$ satisfies $(x_v, y_v) = (1,0)$. This shows that the unsafe edges with respect to $(\bm x, \bm y, \bm z)$ form a bipartite subgraph of $G$. 

	Define $\phi \colon S_= \to S_\times$ as follows. Fix some arbitrary ordering of $V$. For any $(\bm x, \bm y, \bm z) \in S_=$, let $T$ be the lexicographically-first subset of $V$ so that every unsafe edge with respect to $(\bm x, \bm y, \bm z)$ has exactly one endpoint in $T$. Such $T$ exists since the unsafe edges form a bipartite subgraph. Define $\phi(\bm x, \bm y, \bm z) = (\bm x', \bm y', \bm z)$ by setting 
	\[
	(x'_v, y'_v) = \begin{cases}
	 	(y_v, x_v) & \text{if } v \in T, \\
 		(x_v, y_v) & \text{if } v \notin T.
	 \end{cases}
	\]
	In other words, the map $\phi$ swaps $(x_v, y_v)$ for each $v \in T$.
	
	Let us check that the image of $\phi$ lies in $S_\times$, we need to check  that $z_{uv} \leq H(x'_u, y'_v) H(y'_u,x'_v)$ for all $uv \in E$. Only unsafe edges have a chance of violating the inequality. If $uv$ is an unsafe edge, then exactly one of $(x_u, y_u)$ and $(x_v, y_v)$ is swapped by $\phi$, and so $H(x'_u, y'_v) H(y'_u,x'_v) = H(x_u, x_v)H(y_u,y_v)  \ge z_{uv}$. Thus the image of $\phi$ lies in $S_\times$.
	
	To see that $\phi$ is injective, note that given $\phi(\bm x, \bm y, \bm z)$, we can identify the unsafe edges, which are unaffected by swapping, and then recover the lexicographically-first subset $T$ of vertices that contains exactly one vertex from every unsafe edge, and then swap the pair $(x'_v, y'_v)$ for every $v \in T$ to recover $(\bm x, \bm y)$. It is also easy to see that $\phi$ is a measure-preserving map, as we can partition $S_=$ into regions indexed by the set $T$ of swapped vertices. Thus $\phi \colon S_= \to S_\times$ is a measure-preserving injection.
\end{proof}

\section{Colorings} \label{sec:color}

\subsection{Semiproper colorings}

Let us state a generalization of \cref{thm:semiproper} to semiproper list colorings, where every vertex in $G$ has a possibly different set of allowable colors. Recall that ``colors'' are synonymous with vertices of $H = K_q^{\ell\circ}$.

To state the theorem, we will need to set up some notation. Throughout this entire section, we fix a finite set of colors $\Omega$, as well as a subset $\Omega_\circ \subseteq \Omega$ of looped colors. Recall that a \emph{semiproper coloring} of $G$ is an assignment of vertices of $G$ to colors so that no two adjacent vertices of $G$ share the same non-looped color.
The sets of colors $\Omega_\circ \subseteq \Omega$ are fixed throughout, and we omit mentioning them explicitly in the statements below.

For $A,B \subseteq \Omega$ and nonnegative integers $a, b$, define
\[
\cc{A}{B}{a}{b}
\]
to be the number of semiproper colorings of $K_{a,b}$, where the $a$ vertices in the first vertex part of $K_{a,b}$ have their colors chosen from $A$, and the $b$ vertices in the second vertex part of $K_{a,b}$ have their colors chosen from $B$. Observe that $\cc{A}{B}{a}{b}
 = \norm{H|_{A \times B}}_{K_{a,b}}^{ab}$ where $H|_{A\times B}$ is the restriction of the associated partially looped complete graph $H$ (viewed as function $\Omega \times \Omega \to \{0,1\}$) to $A\times B$.

Here is the main theorem of this section. It implies \cref{thm:semiproper} after taking $\Omega_v = \Omega$ for all $v \in V$.

\begin{theorem} \label{thm:semiproper-list}
	Let $G = (V,E)$ be a graph without isolated vertices. Assign a subset of colors $\Omega_v \subseteq \Omega$ to each $v \in V$. Then the total number of semiproper colorings of $G$ where each $v \in V$ is assigned some color from $\Omega_v$ is at most
	\[
		\prod_{uv \in E} \cc{\Omega_u}{\Omega_v}{d_v}{d_u}^{1/(d_ud_v)}.
	\]
\end{theorem}

Here are some conventions about notation that will be maintained throughout this section:
\begin{itemize}
	\item $A \om B := A \setminus (B \setminus \Omega_\circ)$, i.e., remove from $A$ all non-looped colors in $B$. This is a handy operation when we consider what happens to the list of colors at the vertex after we assign colors to its neighbors.
	\item In $A \setminus x$, $A \cup x$, $A \om x$, for $x \in \Omega$, we treat $x$ as a singleton set $\{x\}$.
	\item $\bm x$ and $\bm y$ refer to a vector of colors (colors are elements of $\Omega$), and $x_i$ refers to the $i$-th coordinate of $\bm x$.
	\item Given a vector $\bm x$, we often treat $\bm x$ as a subset of $\Omega$. So $\abs{\bm x}$ is the number of distinct colors appearing in $\bm x$, $y \cup \bm x$ (where $y \in \Omega$) is the union of the elements in $\bm x$ along with $y$, and $A \om \bm x$ is the set of colors left in $A$ after we remove all non-looped colors appearing in $\bm x$.
\end{itemize}

\subsection{Some correlation inequalities for symmetric polynomials}

The main result of this section is the following inequality of symmetric polynomials. We will need it later for our proof of \cref{thm:semiproper-list}.

\begin{proposition} \label{prop:sym}
	Let $\alpha_1, \dots, \alpha_n \ge 0$ be reals, and $k$ a nonnegative integer. Let $\abs{\bm x}$ denote the number of distinct entries in $\bm x$. Set
	\[
	m_\ell := \EEE_{\substack{\bm x \in [n]^k \\ \abs{\bm x} = \ell}} \prod_{i=1}^k \alpha_{x_i}. 
	\]
	Then $m_1 \ge \cdots \ge m_{\min\{n, k\}}$.
\end{proposition}

For example, with $n = 3$ and $k = 4$, we have
\begin{align*}
m_1 &= \tfrac13 (\alpha_1^4 + \alpha_2^4 + \alpha_3^4),
\\
m_2 &= \tfrac{1}{42} (4 \alpha_1^3 \alpha_2 +  4 \alpha_1^3 \alpha_3 + 4 \alpha_2^3 \alpha_1 + 4 \alpha_2^3 \alpha_3 + 4 \alpha_3^3 \alpha_1 + 4 \alpha_3^3 \alpha_2 +  6 \alpha_1^2 \alpha_2^2 + 6 \alpha_1^2 \alpha_3^2+ 6 \alpha_2^2 \alpha_3^2), 
\\
m_3 &= \tfrac13(\alpha_1^2\alpha_2\alpha_3 + \alpha_1\alpha_2^2\alpha_3 + \alpha_1\alpha_2\alpha_3^2).
\end{align*}

For $S \subseteq [n]$ and $\abs{S} \le k$, define
\[
f_{k,S} = \sum_{\substack{\bm x \in S^k \\ \abs{\bm x} = \abs{S}}} \prod_{i=1}^k \alpha_{x_i}.
\]
Here $\abs{\bm x} = \abs{S}$ in the index of the summation simply says that all elements of $S$ appear in $\bm x$. In other words, $f_{k,S}$ is the sum of all monomials whose set of indices is exactly $S$. For example,
\[
f_{5,\{1,2,3\}} = 20 \alpha_1^3\alpha_2\alpha_3 + 20 \alpha_1\alpha_2^3\alpha_3 + 20 \alpha_1\alpha_2\alpha_3^3 + 30 \alpha_1^2\alpha_2^2\alpha_3 + 30 \alpha_1^2\alpha_2\alpha_3^2 + 30  \alpha_1\alpha_2^2\alpha_3^2.
\]
Observe that $f_{k,S}$ satisfies the recursion
\begin{equation}
	\label{eq:f-recur}
f_{k,S} = \sum_{x \in S} \alpha_x (f_{k-1,S} + f_{k-1,S\setminus x}).
\end{equation}

We introduce the following averaging notation. For any polynomial $P$ in the variables $\alpha_1, \dots,$ write $\avg{P} := P / c$ where $c$ is the normalizing constant chosen so that $\avg{P} = 1$ whenever $\alpha_1=\alpha_2 = \cdots = 1$, i.e., $c$ is the sum of all coefficients (if $c = 0$, we set $\avg{P} = 0$). For example, $\avg{\alpha_1 + \alpha_2} = (\alpha_1 + \alpha_2)/2$. For this notation to make sense, we view the $\alpha_1, \dots, \alpha_n$ as formal unassigned variables. When we say that an inequality is true, we mean that it is true for all nonnegative assignments of the $\alpha_i$'s. 
This averaging notation has the convenience that we do not have to keep track of the unimportant normalization factor.

The proof of \cref{prop:sym} proceeds in several steps. 

\begin{lemma} \label{lem:sym-poly-k}
Let $S \subseteq [n]$ and $1 \le \abs{S} < k$. Then
\[
\avg{\sum_{x \in S} \alpha_x f_{k-1,S}} \le \avg{ f_{k,S}}.
\]
\end{lemma}

\begin{proof}
	We apply induction on $\abs{S} + k$. The claimed inequality is an equality when $\abs{S} = 1$ or $k = \abs{S} + 1$. So assume that $\abs{S} > 1$ and $k > \abs{S} + 1$. 
	
	Note that $(\alpha_x : x \in S)$ and $( \alpha_x f_{k-1,S\setminus x} : x \in S)$ are oppositely sorted (meaning, whenever evaluated at nonnegative assignment of the $\alpha_x$'s). Indeed, note that $\alpha_x f_{k-1,S\setminus x} = \alpha_x\alpha_y Q$ where $Q$ is some polynomial with nonnegative coefficients and it does not involve the variable $\alpha_x$, so that if $\alpha_x \le \alpha_y$, then swapping the two variables $\alpha_x$ and $\alpha_y$ cannot increase $Q$. In particular, this sortedness implies, via the rearrangement inequality,
	\[
		\avg{ \sum_{x \in S} \alpha_x^2  f_{k-1,S\setminus x}}
		\le
		\avg{ \sum_{x,y \in S} \alpha_x\alpha_y  f_{k-1,S\setminus x}},
	\]
	and thus, using that $\avg{A} \le \avg{A+B}$ implies $\avg{A+B} \le \avg{B}$, we have
	\begin{equation}
		\label{eq:f-sortedness}
		\avg{ \sum_{x,y \in S} \alpha_x\alpha_y  f_{k-1,S\setminus x}}
		\le
		\avg{ \sum_{x \in S} \sum_{y \in S \setminus x} \alpha_x\alpha_y  f_{k-1,S\setminus x}}.
	\end{equation}
	
	Applying the recursion \eqref{eq:f-recur}, we have
	\[
		\avg{f_{k-1,S}} = \avg{\sum_{x \in S} \alpha_x f_{k-2,S} + \sum_{x \in S} \alpha_x  f_{k-2,S\setminus x}}.
	\]
	On the other hand, by the induction hypothesis with $(k-1, S)$, we have 
	\[
		\avg{\sum_{x \in S} \alpha_x f_{k-2,S}} \le \avg{f_{k-1,S}}.
	\]
	Using that $\avg{A} \le  \avg{A+B}$ implies $\avg{A+B} \le \avg{B}$, we have
	\begin{equation}
		\label{eq:f-k-1-upper}
		\avg{f_{k-1,S}} \le \avg{\sum_{x \in S} \alpha_x f_{k-2,S\setminus x}}.
	\end{equation}
		
	Therefore, we have
	\begin{align*}
		\avg{\sum_{y \in S} \alpha_y f_{k-1,S}}
		& \le \avg{\sum_{x,y \in S} \alpha_x \alpha_y f_{k-2,S\setminus x}} 
		&&\text{\small [by \eqref{eq:f-k-1-upper}]}
		\\
		&\le  \avg{ \sum_{x \in S} \sum_{y \in S \setminus x} \alpha_x\alpha_y  f_{k-2,S\setminus x}}
		&&\text{\small [by \eqref{eq:f-sortedness}]}
		\\
		&\le \avg{ \sum_{x \in S} \alpha_x f_{k-1, S\setminus x}}. && \text{\small[by induction with $(k-1,S\setminus x)$]}
	\end{align*}
	The lemma then follows by using the recursion \eqref{eq:f-recur} and that $\avg{A} \le \avg{B}$ implies $\avg{A} \le \avg{A+B}$.
\end{proof}

\begin{lemma} \label{lem:sym-poly-S}
	Let $S \subseteq [n]$ and $2 \le \abs{S} \le k$. Then
	\[
	\avg{f_{k, S}} \le \avg{\sum_{x \in S} f_{k, S\setminus x}}.
	\]
\end{lemma}

\begin{proof}
	We apply induction on $\abs{S}$. When $\abs{S} = 2$, the lemma follows by noting that $\alpha_1^i\alpha_2^{k-i} + \alpha_1^{k-i}\alpha_2^{i} \le \alpha_1^k + \alpha_2^k$ for all $0 \le i \le k$. Now assume that $\abs{S} > 2$.
	
	Note that $(\alpha_x : x \in S)$ and $(f_{k-1, S\setminus x} : x \in S)$ are oppositely sorted. Indeed, comparing $f_{k-1, S\setminus x}$ with $f_{k-1, S\setminus y}$, we see that $f_{k-1, S\setminus x}$ does not involve $\alpha_x$, and swapping all its $\alpha_y$ to $\alpha_x$ would yield $f_{k-1, S\setminus y}$. Thus, the rearrangement inequality gives 
	\[
	\avg{\sum_{x \in S}\alpha_x f_{k-1,S\setminus x}} \le 	\avg{\sum_{x \in S}\sum_{y \in S\setminus x} \alpha_x f_{k-1,S\setminus y}}.
	\]
	Also, applying the induction hypotheses on $S \setminus x$ for each $x \in S$, we have
	\[
	\avg{ \sum_{x \in S} \alpha_x f_{k-1, S\setminus x}} \le \avg{ \sum_{x \in S}\sum_{y \in S\setminus x} \alpha_x f_{k-1,S\setminus \{x,y\}}}.
	\]
	Using that $\avg{A}\le \avg{B}$ and $\avg{A}\le \avg{C}$ imply $\avg{A}\le \avg{B+C}$, followed by the recursion~\eqref{eq:f-recur}, we obtain
	\[
	\avg{ \sum_{x \in S} \alpha_x f_{k-1, S\setminus x}} \le 
	\avg{\sum_{x \in S}\sum_{y \in S\setminus x} \alpha_x(f_{k-1,S\setminus y} + f_{k-1,S\setminus \{x,y\}})}
	=\avg{\sum_{x \in S} f_{k, S\setminus x}}.
	\]

    \cref{lem:sym-poly-k} and with the recursion \eqref{eq:f-recur} imply (using that $\avg{A}\le \avg{A+B}$ implies $\avg{A+B}\le \avg{B}$)
	\[
	\avg{f_{k,S}} \le  \avg{ \sum_{x \in S} \alpha_x f_{k-1, S\setminus x}}.
	\]
	The Lemma then follows from the above two inequalities.
\end{proof}

\begin{proof}[Proof of \cref{prop:sym}]
Observe that
\[
m_\ell = \EEE_{\substack{S \subseteq [n] \\ \abs{S} = \ell}} \avg{f_{k,S}}.
\]
So \cref{prop:sym} then follows from \cref{lem:sym-poly-S}, as 
\[
m_{\ell-1} = \EEE_{\substack{T \subseteq [n] \\ \abs{T} = \ell-1}} \avg{f_{k, T}}
= \EEE_{\substack{S \subseteq [n] \\ \abs{S} = \ell}} \avg{\sum_{x \in S} f_{k, S\setminus x}}
\ge \EEE_{\substack{S \subseteq [n] \\ \abs{S} = \ell}} \avg{f_{k,S}}
=m_\ell. \qedhere
\]
\end{proof}

\cref{prop:sym} has the following corollary that we will need next.

\begin{corollary} \label{cor:sym}
Let $D$ be a finite set. Let $t \ge 1$ be real. Let $\alpha_x \ge 0$ for each $x \in D$. Let $\abs{\bm x}$ denote the number of distinct elements in $\bm x$. Let $\tau \colon \NN_{\ge 0} \to \RR_{\ge 0}$ be some non-increasing function. Then
\[
\EEE_{\bm x \in D^k} \sqb{ \tau(\abs{\bm x}) } 
\EEE_{\bm x \in D^k} \sqb{\prod_{i=1}^k \alpha_{x_i}}
\le 
\EEE_{\bm x \in D^k} \sqb{  \tau(\abs{\bm x}) \prod_{i=1}^k \alpha_{x_i}}.
\]
\end{corollary}
\begin{proof}
The inequality can be rewritten as
\[
\EEE_{\bm x \in D^k} \sqb{ \tau(\abs{\bm x}) } 
\EEE_{\bm x \in D^k} \sqb{m_{|\bm x|}}
\le \EEE_{\bm x \in D^k} \sqb{  \tau(\abs{\bm x}) m_{|\bm x|}}
\]
which follows from the rearrangement inequality, as $\abs{\bm x}$ and $m_{\abs{\bm x}}$ are oppositely sorted due to \cref{prop:sym}.
\end{proof}

\subsection{Inequalities for semiproper colorings of complete bipartite graphs}

Now we prove some ``local'' inequalities that will be needed in the next section.

\begin{lemma} \label{lem:color-holder}
Let $A, B \subseteq  \Omega$. For any nonnegative integers $k$ and $r \le s \le t$, we have
\[
\cc{A}{B}{k}{s} \le \cc{A}{B}{k}{r}^{\frac{t-s}{t-r}} \cc{A}{B}{k}{t}^{\frac{s-r}{t-r}}.
\]	
\end{lemma}

\begin{proof}
The Lemma follows by H\"older's inequality, after expanding, for each $i \in \{r,s,t\}$,
\[
\cc{A}{B}{k}{i} = \sum_{\bm x \in A^k} \abs{B \om \bm x}^i. \qedhere
\]
\end{proof}

\begin{lemma} \label{lem:color-BCD} 
Let $D \subseteq C \subseteq \Omega$ and $B \subseteq \Omega$, and suppose that $D\subseteq\Omega\setminus\Omega_\circ$. For positive integers $b,c,k$ with $b \ge 2$, and real $t \ge 1$, we have
\[
	\sum_{\bm x \in D^k} \prod_{i=1}^k \cc{B \setminus x_i}{C \setminus x_i}{c-1}{b-1}^{\frac{t}{b-1}}
	\le
	\sum_{\bm x \in D^k} \frac{\abs{C \setminus \bm x}^t} {\abs{C}^{(1 - \frac{k}{c})t}} \prod_{i=1}^k  \cc{B \setminus x_i}{C}{c}{b-1}^{\frac{t (c-1)}{(b-1)c}}.
\]
\end{lemma}

\begin{proof}
We have
\[
\cc{B \setminus x_i}{C \setminus x_i}{c-1}{b-1}
= \sum_{\bm y \in (C\setminus x_i)^{b-1}} \abs{B \om (x_i \cup \bm y)}^{c-1}
= (\abs{C}-1)^{b-1} \EEE_{\bm y \in (C\setminus x_i)^{b-1}} \sqb{\abs{B \om (x_i \cup \bm y)}^{c-1}}.
\]
We used above that $(B\setminus x_i)\om\bm y = B\om (x_i\cup\bm y)$ since $x_i\in D\subseteq\Omega\setminus\Omega_\circ$ is always non-looped. Similarly, 
\[
\cc{B \setminus x_i}{C}{c-1}{b-1}
= \sum_{\bm y \in C^{b-1}} \abs{B \om (x_i \cup \bm y)}^{c-1}
= \abs{C}^{b-1} \EEE_{\bm y \in C^{b-1}} \sqb{\abs{B \om (x_i \cup \bm y)}^{c-1}}.
\]
By linearity of expectations,
\[
\EEE_{\bm x \in C^k} \sqb{ 1-\frac{\abs{\bm x}}{\abs{C}} } = \paren{1 - \frac{1}{\abs{C}}}^k.
\]
It follows that the claimed inequality is equivalent to
\begin{multline} \label{eq:pre-corr}
\bigg(\EEE_{\bm x \in C^k} \sqb{ 1-\frac{\abs{\bm x}}{\abs{C}} }\bigg)^t
\EEE_{\bm x \in D^k} \sqb{
	\prod_{i=1}^k  
		\EEE_{\bm y \in (C\setminus x_i)^{b-1}} 
			\sqb{\abs{B \om (x_i \cup \bm y)}^{c-1}}^{\frac{t}{b-1}}
} 
\\
\le 
\EEE_{\bm x \in D^k} \sqb{ \paren{1-\frac{\abs{\bm x}}{\abs{C}}}^t \prod_{i=1}^k 
\EEE_{\bm y \in C^{b-1}} \sqb{\abs{B \om (x_i \cup \bm y)}^c}^{\frac{t(c-1)}{c(b-1)}}}.
\end{multline}

Applying \cref{cor:sym} with $\alpha_i = \EEE_{\bm y \in C^{b-1}} \sqb{\abs{B \om (x_i \cup \bm y)}^{c-1}}^{\frac{t}{b-1}}$ for each $1 \le i \le k$, we obtain
\begin{multline} \label{eq:corr}
\EEE_{\bm x \in D^k} \sqb{ \paren{1-\frac{\abs{\bm x}}{\abs{C}}}^t }
\EEE_{\bm x \in D^k} \sqb{\prod_{i=1}^k  \EEE_{\bm y \in C^{b-1}} \sqb{\abs{B \om (x_i \cup \bm y)}^{c-1}}^{\frac{t}{b-1}}}
\\
\le 
\EEE_{\bm x \in D^k} \sqb{ \paren{1-\frac{\abs{\bm x}}{\abs{C}}}^t \prod_{i=1}^k \EEE_{\bm y \in C^{b-1}} \sqb{\abs{B \om (x_i \cup \bm y)}^{c-1}}^{\frac{t}{b-1}}}.
\end{multline}

For each $1 \le i \le k$, 
\[
\EEE_{\bm y' \in (C\setminus x_i)^{b-1}} 
			\sqb{\abs{B \om (x_i \cup \bm y')}^{c-1}}
\le 
\EEE_{\bm y \in C^{b-1}} 
			\sqb{\abs{B \om (x_i \cup \bm y)}^{c-1}},
\]
since we can couple  $(\bm y', \bm y)$ so that $x_i \cup \bm y' \supseteq x_i \cup \bm y$ for each $(\bm y', \bm y)$: sample $\bm y$ uniformly from $C^{b-1}$ and obtain $\bm y'$ by replacing every coordinate of $\bm y$ equal to $x_i$ by an independent uniformly random element of $C \setminus x_i$. 

Also, we have
\[
\paren{ \EEE_{\bm x \in D^k} \sqb{1-\frac{\abs{\bm x}}{\abs{C}}}}^t \le \EEE_{\bm x \in D^k} \sqb{ \paren{1-\frac{\abs{\bm x}}{\abs{C}}}^t}.
\]

Combining the above two inequalities, we see that LHS of \cref{eq:corr} upper bounds the LHS of \cref{eq:pre-corr}. Also, the RHS of \cref{eq:corr} lower bounds the RHS of \cref{eq:pre-corr} by the monotonicty of $L^p$ norms. Thus \cref{eq:pre-corr} holds.
\end{proof}

\begin{lemma} \label{lem:color-AC}
Let $A,B,C \subseteq \Omega$. Let $a,b,c$ be positive integers with $\max\{b,c\} \le a$. We have
\[
\sum_{x \in A} \cc{B \om x}{C \om x}{c-1}{b-1}^{\frac{a}{b+c-2}}
	\le 
	\cc{A}{C}{c}{a}^{\frac{b-1}{c(b+c-2)}} \left(\sum_{x \in A} \cc{B\om x}{C}{c}{b-1}^{\frac{a}{c}} \right)^{\frac{c-1}{b+c-2}}.
\]	
\end{lemma}

\begin{proof} When $b = 1$ this is in fact an equality as both sides are equal to $\sum_{x \in A} \abs{B \om x}^a$. From now on assume $b\ge 2$. Raising both sides to exponent $c$ and using $\cc{A}{C}{c}{a} = \sum_{\bm x \in A^c} \abs{C \om \bm x}^a$, the inequality can be rewritten as
\[
\sum_{\bm x \in A^c}\prod_{i=1}^c \cc{B \om x_i}{C \om x_i}{c-1}{b-1}^{\frac{a}{b+c-2}}
\le 
\left(\sum_{\bm x \in A^c} \abs{C \om \bm x}^a \right)^{\frac{b-1}{b+c-2}} 
\left(\sum_{\bm x \in A^c} \prod_{i=1}^c \cc{B \om x_i}{C}{c}{b-1}^{\frac{a}{c}} \right)^{\frac{c-1}{b+c-2}}.
\]
Applying H\"older's inequality to the right-hand side, we see that it suffices to prove
\[
\sum_{\bm x \in A^c}\prod_{i=1}^c \cc{B \om x_i}{C \om x_i}{c-1}{b-1}^{\frac{a}{b+c-2}}
\le 
\sum_{\bm x \in A^c} \abs{C \om \bm x}^{\frac{a(b-1)}{b+c-2}} \prod_{i=1}^c \cc{B \om x_i}{C}{c}{b-1}^{\frac{a(c-1)}{(b+c-2)c}}.
\]

Let $D: = (A \cap C) \setminus \Omega_\circ\subseteq \Omega\setminus \Omega_\circ$. It suffices to show that the above inequality holds with a partial summation where we hold fixed the coordinates of $\bm x$ lying outside $D$ and let the other coordinates range over $D$. In other words, let $K \subseteq [c]$. Fix $x_i \in A\setminus D$ for each $i \notin K$. Writing $\bm x_K = (x_i)_{i \in K}$, it suffices to show that
\begin{equation}\label{eq:c-ABC-reduce}
\sum_{\bm x_K \in D^K}\prod_{i=1}^c \cc{B \om x_i}{C \om x_i}{c-1}{b-1}^{\frac{a}{b+c-2}}
\le 
\sum_{\bm x_K \in D^K} \abs{C \om \bm x}^{\frac{a(b-1)}{b+c-2}} \prod_{i=1}^c \cc{B \om x_i}{C}{c}{b-1}^{\frac{a(c-1)}{(b+c-2)c}}.
\end{equation}

Let $k = |K|$ Applying \cref{lem:color-BCD} with $t = \frac{a(b-1)}{b+c-2} \ge \frac{c(b-1)}{b+c-2} \ge 1$, we obtain
\begin{equation} \label{eq:c-ABC-apply}
\sum_{\bm x_K \in D^K} \prod_{i \in K} \cc{B \setminus x_i}{C \setminus x_i}{c-1}{b-1}^{\frac{a}{b+c-2}}
\le
\sum_{\bm x_K \in D^K} \left(\prod_{i \in K}  \cc{B \setminus x_i}{C}{c}{b-1}^{\frac{a(c-1)}{(b+c-2)c}} \right)
 \frac{\abs{C \setminus \bm x_K}^{\frac{a(b-1)}{b+c-2}}}{\abs{C}^{\frac{a(b-1)(c-k)}{c(b+c-2)}}}.
\end{equation} 

For each $i \notin K$, either $x_i \in \Omega_\circ$ or $x_i \notin C$, so $C \om x_i = C$. Applying \cref{lem:color-holder},
\begin{multline*}
\cc{B \om x_i}{C \om x_i}{c-1}{b-1}
= \cc{B \om x_i}{C }{c-1}{b-1}
\\
\le 
\cc{B \om x_i}{C}{c}{b-1}^{\frac{c-1}{c}} \cc{B \om x_i}{C}{0}{b-1}^{\frac{1}{c}}
= 
\cc{B \om x_i}{C}{c}{b-1}^{\frac{c-1}{c}} \abs{C}^{\frac{b-1}{c}}.
\end{multline*}
By multiplying onto \cref{eq:c-ABC-apply} the above inequality raised to exponent $a/(b+c-2)$ and ranged over all $i\in [c] \setminus K$, we obtain \cref{eq:c-ABC-reduce}. Note that $C \setminus \bm x_K = C \om \bm x$ since $x_i \in D \subseteq C \setminus \Omega_\circ$ for all $i \in K$.
\end{proof}

The following lemma is the ``local'' inequality that the proof of \cref{thm:semiproper-list} will reduce to.

\begin{lemma}[Local inequality for semiproper colorings] \label{lem:color-ABC}
Let $A,B,C \subseteq \Omega$. Let $a,b,c$ be positive integers with $\max\{b,c\} \le a$. We have
\[
	\sum_{x \in A} \cc{B \om x}{C \om x}{c-1}{b-1}^{\frac{a}{b+c-2}}
	\le
	\cc{A}{B}{b}{a}^{\frac{c-1}{b(b+c-2)}} 
	\cc{A}{C}{c}{a}^{\frac{b-1}{c(b+c-2)}} 
	\cc{B}{C}{c}{b}^{\frac{a(b-1)(c-1)}{(b+c-2)bc}}.
\]
\end{lemma}

\begin{proof} 
When $b = 1$ we have an equality: $\sum_{x \in A} \abs{B \om x}^a = \ccs{A}{B}{1}{a}$.
Now let $b \ge 2$. Applying \cref{lem:123} with $\Delta = a$, $\beta=c$, $\gamma = b$, and $f_{12} \colon A \times B \to \{0,1\}$ and $f_{23} \colon B \times C \to \{0,1\}$  the color compatibility functions (i.e., $f(x,y) = 1$ if $x=y \notin \Omega_\circ$ and $f(x,y) = 1$ otherwise; so $\norm{f_{12}}_{K_{s,t}}^{st} = \ccs{A}{B}{s}{t}$ and likewise with $f_{23}$), we obtain
\[
\sum_{x \in A} \cc{B\om x}{C}{c}{b-1}^{\frac{a}{c}} \le \cc{A}{B}{b}{a}^{\frac{1}{b}} \cc{B}{C}{c}{b}^{\frac{a(b-1)}{bc}}.
\]
The Lemma follows by bounding the right-hand side of the inequality in \cref{lem:color-AC} using the above inequality.
\end{proof}

\subsection{Proof of \cref{thm:semiproper-list}} 

Now we are ready to prove \cref{thm:semiproper-list}. We proceed similarly to \cref{sec:pf-BL-main}, with an important new twist, namely that the neighborhood of a vertex is no longer necessarily an independent set. This explains the needs for the more involved inequalities for semiproper colorings seen earlier.  We begin by recalling the notation used in \cref{sec:pf-BL-main} which will also be used in this section.

We apply induction on the number of vertices in $G = (V,E)$. The base case $\abs{E} = 0$ is trivial. Let $\Delta$ be the maximum degree of $G$, and fix a vertex $w$ of degree $\Delta$ in $G$. 

For $k \in \{0, 1, 2, \dots\}\cup\{\infty\}$, let $V_k$ be the set of vertices at distance exactly $k$ from $w$. For $0 \le i \le j \le i+1$, let $E_{ij}$ be the edges with one endpoint in $V_i$ and the other in $V_j$. Let $V_{\ge k} = \bigcup_{i \ge k} V_i$, $E_{\ge k} = \bigcup_{k \le i \le j} E_{ij}$, and $E_{> k} = E_{k,k+1} \cup E_{\ge k+1}$. Note that $V = V_{\ge 0}$ and $E = E_{> 0}$. Note that unlike in the triangle-free setting, now $E_{11}$ may be nonempty.

Let $I_1$ be the set of vertices in $V_1$ whose neighborhood is exactly $\{w\}$.

Although we treat edges as unordered pairs, when we write $vu \in E_{ij}$, we always mean $v \in V_i$ and $u \in V_j$. On the other hand, when we range over $uv \in E_{ii}$, we do not count $uv$ and $vu$ separately.

We use the variable $x_v$ to denote an element of $\Omega_v$. For any $S \subseteq V$, write $\Omega_S := \prod_{v \in S} \Omega_v$ and $d \bm x_S := \prod_{v \in S} dx_v$.

We associate to each edge $uv$ the function $f_{uv} \colon \Omega_u \times \Omega_v \to \{0,1\}$ encoding validity of color assignments:
\[
f_{uv} (x_u, x_v) = \begin{cases} 0 & \text{if } x_u = x_v \notin \Omega_\circ, \\
1 & \text{otherwise.}
\end{cases}
\]

Define
\[
f_{vu}^{x_w}(x_v,x_u) := f_{wv}(x_w,x_v)^{1/(d_v - 1)} f_{vu}(x_v, x_u), \quad \text{for } vu \in E_{12} \cup E_{12},
\]
and
\[
f_{vu}^{x_w,x_w}(x_v,x_u) := f_{wv}(x_w,x_v)^{1/(d_v - 1)} f_{wv}(x_w,x_u)^{1/(d_u - 1)} f_{vu}(x_v, x_u), \quad \text{for } vu \in E_{11}.
\]

\smallskip 

We have
\begin{align*}
	&\int_{\Omega_V} \prod_{vu\in E} f_{vu}(x_v,x_u) \, d\bm x_V
	\\ 
	&= \int_{\Omega_{V\setminus I_1}} \prod_{v \in I_1} \left(\int_{\Omega_v} f_{wv}(x_w,x_v) \, dx_v\right) 
	\prod_{vu\in E_{11}} f_{vu}^{x_w,x_w}(x_v,x_u) 
	\prod_{vu\in E_{12}} f_{vu}^{x_w}(x_v,x_u) 
	 \prod_{vu\in E_{\ge 2}} f_{vu}(x_v,x_u) \,
	 d\bm x_{V\setminus I_1}
	\\
	&\le \int_{\Omega_w} \prod_{v \in I_1} 
	\norm{f_{wv}(x_w,\cdot)}_1 
	\prod_{vu\in E_{11}} \norm{f_{vu}^{x_w,x_w}}_{K_{d_u-1,d_v-1}} 
	\prod_{vu\in E_{12}} \norm{f_{vu}^{x_w}}_{K_{d_u,d_v-1}} \prod_{vu\in E_{\ge 2}} \norm{f_{vu}}_{K_{d_u,d_v}} \, d x_w,
\end{align*}
where in the last step we applied the induction hypothesis to $G-w$. 
It remains to prove the bound
\begin{multline*}
\int_{\Omega_w} \prod_{v \in I_1} 
	\norm{f_{wv}(x_w,\cdot)}_1 
	\prod_{vu\in E_{11}} \norm{f_{vu}^{x_w,x_w}}_{K_{d_u-1,d_v-1}} 
	\prod_{vu\in E_{12}} \norm{f_{vu}^{x_w}}_{K_{d_u,d_v-1}} \prod_{vu\in E_{\ge 2}} \norm{f_{vu}}_{K_{d_u,d_v}} \, d x_w
\\
\le \prod_{vu\in E} \|f_{vu}\|_{K_{d_u,d_v}}.	
\end{multline*}
\textbf{First localization.}
Observing that the factor $\norm{f_{vu}}_{K_{d_u,d_v}}$ appears on both sides whenever $vu \in E_{\ge 2}$, we see that it suffices to prove
\begin{multline*}
\int_{\Omega_w} \prod_{v \in I_1} 
	\norm{f_{wv}(x_w,\cdot)}_1 
	\prod_{vu\in E_{11}} \norm{f_{vu}^{x_w,x_w}}_{K_{d_u-1,d_v-1}} 
	\prod_{vu\in E_{12}} \norm{f_{vu}^{x_w}}_{K_{d_u,d_v-1}} 
	\, d x_w
	\\
	\le \prod_{vu\in E_{01} \cup E_{11} \cup E_{12}} \|f_{vu}\|_{K_{d_u,d_v}}.
\end{multline*}
By distributing the $E_{01}$ factors on the RHS to $E_{11} \cup E_{12}$, we can rewrite the above inequality as
\begin{multline}
	\label{eq:c-after-trim}
	\int_{\Omega_w} \prod_{v \in I_1} 
	\norm{f_{wv}(x_w,\cdot)}_1 
	\prod_{vu\in E_{11}} \norm{f_{vu}^{x_w,x_w}}_{K_{d_u-1,d_v-1}} 
	\prod_{vu\in E_{12}} \norm{f_{vu}^{x_w}}_{K_{d_u,d_v-1}} 
	\, d x_w
	\\
	\le 
	\prod_{v \in I_1} \|f_{wv}\|_{K_{1,\Delta}} 
	\prod_{vu \in E_{11}} \left( \|f_{wv}\|_{K_{d_v,\Delta}}^{1/(d_v-1)} \|f_{wv}\|_{K_{d_u,\Delta}}^{1/(d_u-1)} \|f_{vu}\|_{K_{d_u,d_v}}\right)
	\\
	\cdot \prod_{vu \in E_{12}} \left( \|f_{wv}\|_{K_{d_v,\Delta}}^{1/(d_v-1)} \|f_{vu}\|_{K_{d_u,d_v}}\right).
\end{multline}
\textbf{Second localization.}
Applying H\"older's inequality with exponents given by the summands of 
\[
\sum_{v \in I_1} \frac{1}{\Delta} + \sum_{vu \in E_{11}} \left(\frac{1}{\Delta(d_v-1)} + \frac{1}{\Delta(d_u-1)}\right) +  \sum_{vu \in E_{12}} \frac{1}{\Delta(d_v-1)}= 1,
\]
we upper bound the left-hand side of \cref{eq:c-after-trim} by
\begin{multline}
\int_{\Omega_w} \prod_{v \in I_1} \norm{f_{wv}(x_w,\cdot)}_1 
\prod_{vu\in E_{11}} \norm{f_{vu}^{x_w,x_w}}_{K_{d_u-1,d_v-1}} 
 \prod_{vu\in E_{12}} \norm{f_{vu}^{x_w}}_{K_{d_u,d_v-1}} \, d x_w\
\\
\le 
\prod_{v \in I_1} \left( \int_{\Omega_w} \norm{f_{wv}(x_w,\cdot)}_1^\Delta \, d x_w\right)^{\frac{1}{\Delta}}
\prod_{vu\in E_{11}} \left( \int_{\Omega_w} \norm{f_{vu}^{x_w,x_w}}_{K_{d_u-1,d_v-1}}^{\frac{\Delta(d_v-1)(d_u-1)}{d_v+d_u-2}
} \, d x_w\right)^{\frac{d_v+d_u-2}
{\Delta(d_v-1)(d_u-1)}} 
\\
\cdot \prod_{vu\in E_{12}} \left( \int_{\Omega_w} \norm{f_{vu}^{x_w}}_{K_{d_u,d_v-1}}^{\Delta(d_v-1)} \, d x_w\right)^{\frac{1}
{\Delta(d_v-1)}}. \label{eq:c-after-trim-LHS-bd}
\end{multline}

We have
\[
\int_{\Omega_w} \norm{f_{wv}(x_w,\cdot)}_1^\Delta \, d x_w = \norm{f_{wv}}_{K_{1,\Delta}}^\Delta.
\]

For every $vu \in E_{12}$, applying \cref{lem:123} yields
\[
	\int_{\Omega_w} \norm{f_{vu}^{x_w}}_{K_{d_u,d_v-1}}^{\Delta(d_v-1)} \, d x_w  
	\le \|f_{wv}\|_{K_{d_v,\Delta}}^\Delta  \|f_{vu}\|_{K_{d_u,d_v}}^{\Delta(d_v-1)}.
\]

For every $vu \in E_{11}$, applying \cref{lem:color-ABC} with  $a = \Delta$, $b = d_v$, $c = d_u$, noting that $\norm{f_{uv}}_{K_{s,t}}^{st} = \ccs{\Omega_u}{\Omega_v}{s}{t}$ (and likewise with other edges), we obtain
\begin{equation}
	\label{eq:c-pf-uv}
\int_{\Omega_w} \norm{f_{vu}^{x_w,x_w}}_{K_{d_u-1,d_v-1}}^{\frac{\Delta(d_v-1)(d_u-1)}{d_v+d_u-2}
} \, d x_w
\le 
\norm{f_{wv}}_{K_{d_v,\Delta}}^{\frac{\Delta(d_u-1)}{d_v+d_u-2}} 
\norm{f_{wu}}_{K_{d_u,\Delta}}^{\frac{\Delta(d_v-1)}{d_v+d_u-2}} 
\norm{f_{vu}}_{K_{d_u,d_v}}^{\frac{\Delta(d_v-1)(d_u-1)}{d_v+d_u-2}}.
\end{equation}
Applying the above three inequalities to upper bound the right-hand side of \cref{eq:c-after-trim-LHS-bd}, we obtain \cref{eq:c-after-trim}, and thereby complete the proof of \cref{thm:semiproper-list}. \qed  

\section{Clique maximizers for positive semidefinite models} \label{sec:clique}

In this section we prove \cref{thm:clique-max} that for every positive semidefinite (i.e., ferromagnetic) $H \colon \Omega \times  \Omega \to \RR_{\ge 0}$, one has
\[
	\hom(G, H) \le \prod_{v \in V(G)} \hom(K_{d_v+1}, H)^{1/(d_v+1)}.
\]

We introduce some notation that allows us to alter the vertex weights of $H$ separately for each vertex of $G$. It can be thought of as a weighted version of list coloring, where each vertex of $G$ has a different vector of weights on the set of ``colors'' (the vertices of $H$).
Given graph $G = (V,E)$ and weighted graph $H \colon \Omega \times \Omega \to \RR_{\ge 0}$, and vector $\bm \lambda = (\lambda_v)_{v \in V}$ whose coordinates $\lambda_v \colon \Omega \to \RR_{\ge 0}$ are measurable functions, write
\[
\hom_{\bm \lambda} (G,H) := \int_{\Omega^{V}} \prod_{uv \in E} H(x_u,x_v) \, \prod_{v \in V} \lambda_v(x_v) dx_v.
\]
For a single $\lambda \colon \Omega \to \RR_{\ge 0}$, we write $\hom_{\lambda} (G,H)$ to mean $\hom_{\bm \lambda} (G,H)$ with $\lambda_v = \lambda$ for all $v \in G$. 
Observe that when $H = K_q$ and every $\lambda_v$ is the indicator function on a subset of $V(H)$, the above quantity is precisely the usual list coloring.

The following theorem generalizes \cref{thm:clique-max} by taking $\lambda_v = 1$ for all $v \in V$.

\begin{theorem} \label{thm:clique-max-list}
	Let $G = (V,E)$ be a graph and $H \colon \Omega \times \Omega \to \RR_{\ge 0}$ a positive semidefinite weighted graph. Let $\bm \lambda = (\lambda_v)_{v \in V}$ where $\lambda_v \colon \Omega \to \RR_{\ge 0}$ is a measurable function. Then
	\[
	\hom_{\bm \lambda} (G,H) \le \prod_{v \in V} \hom_{\lambda_v} (K_{d_v+1}, H)^{1/(d_v+1)}.
	\]
\end{theorem}

Note that the positive semidefiniteness hypothesis is optimal, since for $G = K_2$ the above inequality is just the Cauchy--Schwarz inequality for the bilinear form in $\RR^\Omega$ given by $H$.

Our proof of \cref{thm:clique-max-list} proceeds by induction on the maximum degree of $G$.
We begin with a lemma that is essentially an application of the Cauchy--Schwarz inequality.
Let $G^\bullet$ be the graph obtained from $G$ by adding a new vertex adjacent to all other vertices. Let $G^{\bullet\bullet}$ be the graph obtained from $G$ by adding two new vertices adjacent to all vertices of $G$ but not to each other.

\begin{lemma} \label{lem:clique-CS}
	Let $G = (V,E)$ be a graph and $H \colon \Omega \times \Omega \to \RR_{\ge 0}$ a weighted graph. Let $\bm \lambda = (\lambda_v)_{v \in V(G^{\bullet\bullet})}$, $\bm \mu = (\mu_v)_{v \in V(G)}$, and $\bm \nu = (\nu_v)_{v \in V(G^{\bullet})}$ be vectors of measurable functions $\Omega \to \RR_{\ge 0}$ such that for every $v \in V(G)$, one has $\lambda_v \mu_v = \nu_v^2$. Furthermore, assume that the entries of $\bm \lambda$ and $\bm \nu$ associated to the new vertices (i.e., the vertices not in $V$) are all identical. Then
	\[
	\hom_{\bm \lambda}(G^{\bullet\bullet}, H) \hom_{\bm \mu}(G, H) \ge \hom_{\bm \nu}(G^\bullet, H)^2.
	\]
\end{lemma}

\begin{proof}
	For $\bm x_V = (x_v)_{v \in V} \in \Omega^V$, we write
	\[
	H^G(\bm x_V) := \prod_{uv \in E} H(x_u, x_v),
	\]
	and
	\[
	\bm\lambda(\bm x_V) = \prod_{v \in V} \lambda_v(x_v), \quad
	\bm\mu(\bm x_V) = \prod_{v \in V} \mu_v(x_v), \quad \text{and}\quad 
	\bm\nu(\bm x_V) = \prod_{v \in V} \nu_v(x_v).
	\]
	Recall that the entries of $\bm \lambda$ and $\bm \nu$ associated to the new vertices are all identical, which we call $\nu_\bullet \colon \Omega \to \RR_{\ge 0}$. We have
	\begin{align*}
	\hom_{\bm \lambda}(G^{\bullet\bullet}, H) &= \int_{\Omega^V} \paren{\int_\Omega \prod_{v \in V} H(x_v, y) \,  \nu_\bullet(y)dy}^2 H^G (\bm x_V) \bm\lambda(\bm x_V) \, d\bm x_V,
	\\
	\hom_{\bm \mu}(G, H) &= \int_{\Omega^V}  H^G (\bm x_V) \bm\mu(\bm x_V) \, d\bm x_V, \quad \text{and}
	\\
	\hom_{\bm \nu}(G^{\bullet}, H) &= \int_{\Omega^V} \paren{\int_\Omega \prod_{v \in V} H(x_v, y) \,  \nu_\bullet(y)dy} H^G (\bm x_V) \bm\nu(\bm x_V) \, d\bm x_V.
	\end{align*}
	The Lemma then follows from the Cauchy--Schwarz inequality applied with respect to $d \bm x_V$, noting that $\bm\lambda(\bm x_V)\bm\mu(\bm x_V) = \bm\nu(\bm x_V)^2$.
\end{proof}

Set 
\[
h_a(\lambda) := \hom_\lambda(K_a, H).
\]
In particular, $h_0(\lambda) = 1$.

\begin{lemma} \label{lem:h-log-convex}
Suppose that \cref{thm:clique-max-list} holds for all $G$ with maximum degree less than $\Delta$. 
Let $2 \le t < \Delta$ be a positive integer. Let $\lambda, \mu, \nu \colon \Omega \to \RR_{\ge 0}$ be measurable functions satisfying $\lambda \mu = \nu^2$ pointwise. Then
\[
h_{t+1}(\lambda)^{1/(t+1)}h_{t-1}(\mu)^{1/(t-1)}
\ge h_t(\nu)^{2/t}.
\]
\end{lemma}

\begin{proof}
	Define $\bm \lambda' = (\lambda'_v)_{v \in V(K_{t-1}^{\bullet\bullet})}$ by setting $\lambda'_v = \lambda$ for each vertex $v$ in the original $K_{t-1}$ of $K_{t-1}^{\bullet\bullet}$ and $\lambda'_v = \nu$ for the two other vertices $v$ (each of degree $t-1$). By \cref{lem:clique-CS},
	\[
	\hom_{\bm \lambda'}(K_{t-1}^{\bullet\bullet}, H) \hom_{\mu}(K_{t-1}, H) \ge \hom_{\nu}(K_{t-1}^{\bullet}, H)^2.
	\]
	Since the maximum degree of $K_{t-1}^{\bullet\bullet}$ is $t < \Delta$, using the hypothesis of the Lemma to apply \cref{thm:clique-max-list} to $K_{t-1}^{\bullet\bullet}$, we have
	\[
	\hom_{\bm \lambda'}(K_{t-1}^{\bullet\bullet}, H) \le h_{t+1}(\lambda)^{\frac{t-1}{t+1}} h_{t}(\nu)^{\frac{2}{t}},
	\]
	so that the previous inequality implies that
	\[
		h_{t+1}(\lambda)^{\frac{t-1}{t+1}} h_t(\nu)^{\frac{2}{t}} h_{t-1}(\mu) \ge h_{t}(\nu)^2.
	\] 
	The Lemma follows after rearranging.
\end{proof}

\begin{remark}
The hypothesis $t < \Delta$ in \cref{lem:h-log-convex} is important for applying the induction hypothesis.
\end{remark}

\begin{lemma} \label{lem:F-log-conv}
	Let $H \colon \Omega \to \RR_{\ge 0}$ be positive semidefinite.
	Let $a$ be a positive integer. Let $\mu, \nu \colon \Omega \to \RR_{\ge 0}$ be measurable functions. For each $0 \le i \le a$, let
	\[
	F_i = \hom_{\bm \lambda}(K_a, H),
	\quad
	\text{where }
	\bm\lambda = (\underbrace{\mu, \dots, \mu}_{i \text{ times}}, \underbrace{\nu, \dots, \nu}_{a-i \text{ times}}).
	\] 
	Then $F_0, \dots, F_a$ is log-convex ($F_s F_{s+2} \ge F_{s+1}^2$ for each $0 \le s \le a-2$), and consequently, $F_0^{a-1}  F_{a} \ge F_1^a$.
\end{lemma}

\begin{proof}
    Let $0 \le s \le a-1$. Write $T = [a]\setminus\{s+1,s+2\}$, $\bm x_T = (x_i)_{i \in T} \in \Omega^{a-2}$,
	\[\ell_{1}(x; \bm x_T) = \nu(x) \prod_{i\in T}H(x,x_i),\]
	\[\ell_{2}(x; \bm x_T) = \mu(x) \prod_{i\in T}H(x,x_i),\]
	and 
	\[g = g(\bm x_T)= \prod_{i,j \in T: i \ne j} H(x_i, x_j) \prod_{i=1}^s \mu(x_i)  \prod_{i=s+3}^{a} \nu(x_i).\]
	For $i,j\in\{1,2\}$, write
	\[
	L_{ij} = L_{ij}(\bm x_T) = \int_{\Omega^{2}}H(x,y)\ell_i(x; \bm x_T) \ell_j(y; \bm x_T)\, dxdy.
	\]
	By the positive semidefiniteness of $H$, we have, for every $\bm x_T \in \Omega^{a-2}$,
	\[
	    L_{11} L_{22} \ge L_{12}^2.
	\]
	We can write 
	\[
		F_s = \int_{\Omega^{a-2}} g L_{11}\, d\bm x_T, \quad
		F_{s+1} = \int_{\Omega^{a-2}} g L_{12}\, d\bm x_T, \quad \text{and} \quad 
		F_{s+2} = \int_{\Omega^{a-2}} g L_{22}\, d\bm x_T.
	\]
    By the Cauchy--Schwarz inequality, we have
\[
F_sF_{s+2}
\ge \paren{ \int_{\Omega^{a-2}}g \sqrt{L_{11}L_{22}} \, \bm x_T }^2 
\ge \paren{ \int_{\Omega^{a-2}}g L_{12} \, \bm x_T }^2  = F_{s+1}^2. \qedhere
\]
\end{proof}

\begin{lemma} \label{lem:M-log-conv}
	Let $b \le a \le  \Delta$ be positive integers.  
	Suppose that \cref{thm:clique-max-list} holds for all graphs $G$ with maximum degree less than $\Delta$. 
	Let $H \colon \Omega \times \Omega \to \RR_{\ge 0}$ be positive semidefinite.
	Let $\lambda, \mu \colon \Omega \to \RR_{\ge 0}$.
	Define $\eta \colon \Omega \to \RR_{\ge 0}$ by 
	\[
	\eta(x) := h_b(\mu H(x, \cdot))^{1/b}
	\]
	(here $\mu H(x, \cdot)$ is the pointwise product of two functions $\Omega \to \RR_{\ge 0}$). For each $0 \le s \le a+1$, set
	\[
	M_s := h_s(\lambda \eta^{a+1-s}) h_{b+1}(\mu)^{\frac{s(s-1)}{b+1}}.
	\]
	Then
	\begin{enumerate}
	\item[(a)] $M_{b+1} \ge M_b M_1$ provided that $b < \Delta$;
	\item[(b)] $M_s M_{s+2} \ge M_{s+1}^2$ for all $b \le s \le a-1$ provided that $b < \Delta$;
	\item[(c)] $M_{a+1} \ge M_1^{a+1}$.
	\end{enumerate}
\end{lemma}

\begin{proof}
	(a) Rewriting the desired inequality, we wish to prove
	\[
	h_{b+1}(\lambda \eta^{a-b}) h_{b+1}(\mu)^b 
	\ge 
	h_b(\lambda \eta^{a+1-b}) h_{b+1}(\mu)^{\frac{b(b-1)}{b+1}}
	h_1(\lambda \eta^{a}).
	\]
	Let $\lambda' = \lambda \eta^{a-b}$. The desired inequality can be rewritten as
	\begin{equation} \label{eq:cml-lambda'}
	h_{b+1}(\lambda') h_{b+1}(\mu)^{\frac{2b}{b+1}} 
	\ge 
	h_b(\lambda' \eta) h_1(\lambda' \eta^{b}).
	\end{equation}
	Note that
	\begin{equation}
		\label{eq:h1}
	 h_1(\lambda' \eta^{b}) = \int_\Omega \lambda'(x) h_b(\mu H(x,\cdot)) \, dx = \hom_{(\lambda', \mu, \mu, \dotsc, \mu)} (K_{b+1}, H).
	\end{equation}
	Applying \cref{lem:F-log-conv}, we have
	\begin{equation}
		\label{eq:h_b+1_lambda_mu}
	h_{b+1}(\lambda')^{\frac{1}{b+1}} h_{b+1}(\mu)^{\frac{b}{b+1}} \ge h_1(\lambda' \eta^{b}).
	\end{equation}
	Thus \eqref{eq:cml-lambda'} reduces to (after eliminating $h_{b+1}(\mu)^{\frac{2b}{b+1}}$)
	\begin{equation}
		\label{eq:h_b+1_1}
	h_{b+1}(\lambda')^{\frac{b-1}{b+1}}  h_1(\lambda' \eta^{b})
	\ge 
	h_b(\lambda' \eta).
	\end{equation}
	Let $N_s := h_{b+1-s}(\lambda' \eta^s)^{\frac{1}{b+1-s}}$. \cref{lem:h-log-convex} implies that $N_s$ is log-convex, i.e., $N_s N_{s+2} \ge N_{s+1}^{2}$ for all $0 \le s \le b-2$ (the hypothesis of \cref{lem:h-log-convex} is satisfied since $b < \Delta$), and consequently, $N_0^{b-1} N_b \ge N_1^b$, which proves \eqref{eq:h_b+1_1}.
	
	(b) The desired inequality $M_s M_{s+2} \ge M_{s+1}^2$ is equivalent, upon expanding and simplying, to
	\begin{equation}
		 \label{eq:pf-M-log-conv-h}
	\sqrt{h_s(\nu_1) h_{s+2}(\nu_2)} h_{b+1}(\mu)^{\frac{1}{b+1}}	\ge 	h_{s+1}(\sqrt{\nu_1\nu_2}).
	\end{equation}
	where $\nu_1 = \lambda \eta^{a+1-s}$ and $\nu_2 =  \lambda \eta^{a-1-s}$.
	
	In this proof, we write, for any $\bm x = (x_1, \dots, x_\ell) \in \Omega^\ell$,
	\[
	H(\bm x; y) = \prod_{i=1}^\ell H(x_i, y),
	\quad
	H(\bm x) = \prod_{1 \le i < j \le \ell} H(x_i, x_j),
	\quad 
	\nu(\bm x) = \prod_{i=1}^\ell \nu(x_i),
	\quad
	\text{and} \quad 
	d\bm x = \prod_{i=1}^\ell d x_i.
	\]
	For any measure $\mu$ on $\Omega$, by expanding $h_s(\nu) = \hom_\nu(K_s, H)$ via fixing all but $b-1$ vertices of $K_s$,
	\[
	h_s(\nu) = \int_{\Omega^{s-b+1}} h_{b-1}(\nu H(\bm x; \cdot )) H(\bm x) \nu(\bm x) d \bm x.
	\]
	By the Cauchy--Schwarz inequality,
	\[
	\sqrt{h_s(\nu_1)h_{s+2}(\nu_2)}
	\ge \int_{\Omega^{s-b+1}} 
		\sqrt{
			h_{b-1}(\nu_1 H(\bm x;\cdot))
			h_{b+1}(\nu_2 H(\bm x;\cdot))
			(\nu_1\nu_2)(\bm x)}
			H(\bm x) 
			d \bm x.
	\]
    Also,
    \[
	h_{s+1}(\sqrt{\nu_1\nu_2})
	= \int_{\Omega^{s-b+1}} 
			h_{b}(\sqrt{\nu_1\nu_2} H(\bm x;\cdot))
			\sqrt{(\nu_1\nu_2)(\bm x)}
			H(\bm x) 
			d \bm x.
	\]
	Thus, to establish \cref{eq:pf-M-log-conv-h}, it suffices to show that for every fixed $\bm x\in \Omega^{s-b+1}$,
	\[
	\sqrt{
			h_{b-1}(\nu_1 H(\bm x;\cdot))
			h_{b+1}(\nu_2 H(\bm x;\cdot))
			}
		h_{b+1}(\mu)^{\frac{1}{b+1}} \ge h_b(\sqrt{\nu_1\nu_2}H(\bm x; \cdot)),
	\]
	as integrating it against $\sqrt{(\nu_1\nu_2)(\bm x)} H(\bm x)$ over $\bm x \in \Omega^{s-b+1}$ would yield \cref{eq:pf-M-log-conv-h}.
	
	Writing $\lambda' = \nu_2 H(\bm x; \cdot)$, and recall that $\nu_2 = \nu_1\eta^2$, it remains to prove that 
	\[
	h_{b-1}(\lambda' \eta^2)
	h_{b+1}(\lambda' ) h_{b+1}(\mu)^{\frac{2}{b+1}}
 	\ge 	h_b(\lambda' \eta )^2,
	\]
	which follows from multiplying together the following three inequalities (our earlier proofs establish their validity for all measures $\lambda'$ on $\Omega$):
	\begin{align*}
		h_{b-1}(\lambda' \eta^2)^{\frac{1}{2}} h_{b+1}(\lambda')^{\frac{b-1}{2(b+1)}} &\ge h_b(\lambda'\eta)^{\frac{b-1}{b}}, \text{ and} && \text{[by \cref{lem:h-log-convex}]}
		\\
		h_{b+1}(\lambda')^{\frac{1}{b(b+1)}} h_{b+1}(\mu)^{\frac{1}{b+1}} &\ge h_1(\lambda'\eta^b)^{\frac{1}{b}},  && \text{[by \eqref{eq:h_b+1_lambda_mu}]}
		\\
		h_{b+1}(\lambda')^{\frac{b-1}{b(b+1)}} h_1(\lambda' \eta^b)^{\frac{1}{b}} &\ge h_b(\lambda'\eta)^{\frac{1}{b}}.  && \text{[by \eqref{eq:h_b+1_1}]}
	\end{align*}
	
	(c) We apply induction on $a-b$. When $a = b$, we have
	\[
	M_{b+1} = h_{b+1}(\lambda) h_{b+1}(\mu)^b
	\quad \text{and} \quad 
	M_1 = h_1(\lambda \eta^b).
	\]
	We have $h_{b+1}(\lambda) h_{b+1}(\mu)^b \ge h_1(\lambda \eta^b)^{b+1}$ by \cref{lem:F-log-conv} (noting \eqref{eq:h_b+1_lambda_mu}), and thus $M_{b+1} \ge M_1^b$, as claimed.
	
	Now assume $b < a$. By the induction hypothesis, we have $M_a \ge M_1^a$. From (a) and (b) above (note that $b < \Delta$ now), we have $M_{a+1}/M_a \ge M_a/M_{a-1} \ge \cdots \ge M_{b+1}/M_b \ge M_1$, and thus $M_{a+1} \ge M_1 M_a \ge M_1^{a+1}$, as claimed.
\end{proof}

Now we are ready to prove \cref{thm:clique-max-list}, which, as a reminder, says that for $G = (V,E)$, $\bm \lambda = (\lambda_v : \Omega \to \RR_{\ge 0})_{v \in V}$, and positive semidefinite $H \colon \Omega \times \Omega \to \RR_{\ge 0}$, one has
\begin{equation}
	\label{eq:clique-max-list-h}
	\hom_{\bm \lambda} (G,H) \le \prod_{v \in V} h_{d_v+1}(\lambda_v)^{1/(d_v+1)}.
\end{equation}
	
\begin{proof}[Proof of \cref{thm:clique-max-list}]
	We apply induction first on $\Delta$, an upper bound on the maximum degree of $G$, and then on the number of vertices of $G$. The base case, for each $\Delta$, is when $G$ has no vertices, in which case the statement is trivial. The only non-trivial case is if $G$ contains a vertex $w$ of degree exactly $\Delta$. As earlier, let $V_i$, $i \in \{0, 1, \dots \} \cup \{\infty\}$, denote the set of vertices at distance exactly $i$ from $w$.
	
	First choosing the color on $w$, we obtain
	\[
	\hom_{\bm \lambda} (G, H) = \int_{\Omega} \lambda_w(x_w) \hom_{\bm \mu^{x_w}} (G-w, H) \, dx_w,
	\]
	where $\bm \mu^{x_w} = (\mu_v)_{v \in V(G -w)}$ is defined by $\mu_v = \lambda_v H(x_w, \cdot)$ for $v \in V_1$, and $\mu_v = \lambda_v$ for all $v \in V_{\ge 2}$. Now, by applying the induction hypothesis on $G-w$ to upper bound the integrand, we have,
	\[
	\hom_{\bm \lambda} (G, H) \le 
	\int_{\Omega} \lambda_w(x_w) \prod_{v \in V_1} h_{d_v}(\lambda_v H(x_w, \cdot))^{1/d_v} \prod_{v \in V_{\ge 2}} h_{d_v + 1}(\lambda_v)^{1/(d_v+1)} \, dx_w.
	\]
	Comparing with the right-hand side of \eqref{eq:clique-max-list-h}, we see that it suffices to prove that
	\[
	 \int_{\Omega} \lambda_w(x_w) \prod_{v \in V_1} h_{d_v}(\lambda_v H(x_w, \cdot))^{1/d_v} \, dx_w
	\le 
	 h_{d_w+1}(\lambda_w)^{1/(d_w+1)} \prod_{v \in V_1} h_{d_v+1}(\lambda_v)^{1/(d_v+1)}.
	\]
	Applying H\"older's inequality to the left-hand side (noting that $\abs{V_1} = d_w$), we have
	\begin{align*}
	 \int_{\Omega} \lambda_w(x_w) \prod_{v \in V_1} h_{d_v}(\lambda_v H(x_w, \cdot))^{1/d_v} \, dx_w
	&\le 
	  \prod_{v \in V_1}  \paren{\int_{\Omega} \lambda_w(x_w) h_{d_v}(\lambda_v H(x_w, \cdot))^{d_w/d_v} \, dx_w}^{1/d_w}
	 \\
	 &=
	 \prod_{v \in V_1} h_1(\lambda_w \eta_v^{d_w})^{1/d_w},
	\end{align*}
	where $\eta_v(x) = h_{d_v}(\lambda_vH(x, \cdot))^{1/d_v}$. Thus it suffices to prove that, for each $v \in V_1$,
	\[
	h_1(\lambda_w \eta_v^{d_w})
	\le 
	h_{d_w+1}(\lambda_w)^{1/(d_w+1)} h_{d_v+1}(\lambda_v)^{d_w/(d_v+1)}.
	\]
	But this is exactly \cref{lem:M-log-conv}(c) with $a = d_w$, $b= d_v$, $\lambda = \lambda_w$, $\mu = \lambda_v$, $\eta = \eta_v$. Note that this is a valid application of \cref{lem:M-log-conv}(c) as $d_v \le d_w \le \Delta$, and \cref{thm:clique-max-list} is assumed to hold for all graphs with maximum degree less than $\Delta$ by inductive hypothesis.
\end{proof}

\textbf{Acknowledgments.} We thank Suvrit Sra for suggesting the proof of \cref{lem:mixed-norm}. We would also like to thank the anonymous referees for detailed comments which helped to improve the paper.


\end{document}